\documentclass[11pt]{amsart}

\newif\ifcomments
\commentsfalse

\usepackage[hmarginratio={1:1}]{geometry}
\usepackage[shortlabels]{enumitem}
\usepackage{mathrsfs}
\usepackage[hidelinks]{hyperref}
\usepackage[style=alphabetic,
doi=false,
isbn=false,
maxnames=99,
maxalphanames=99,
backref=false]{biblatex}

\usepackage[dvipsnames]{xcolor}
\usepackage{amssymb}
\usepackage{soulutf8}
\usepackage{tikz-cd}
\usepackage{calc}
\usepackage{mathtools}

\ifcomments
\usepackage{todonotes}
\fi

\newcommand\abs[2][]{#1\lvert{#2}#1\rvert}
\newcommand\AI{\operatorname{AI}}
\newcommand\alttext{\texorpdfstring}
\newcommand\Ar{\operatorname{Ar}}
\newcommand\bA{\mathbf{A}}
\newcommand\basic{\mathsf{basic}}
\newcommand\bB{\mathbf{B}}
\newcommand\Bo{\mathsf{Bo}}
\newcommand\br[1]{\langle{}#1\rangle{}}
\newcommand\cB{\mathcal{B}}
\newcommand\cBo{\mathsf{cBo}}
\newcommand\cC{\mathcal{C}}
\newcommand\cD{\mathcal{D}}
\newcommand\cE{\mathcal{E}}
\newcommand\cL{\mathcal{L}}
\newcommand\cl{\mathsf{cl}}
\newcommand\cm{\mathrm{c}}
\newcommand\code{\mathsf{code}}
\newcommand\con{\triangleleft}
\newcommand\cons{\mathsf{cons}}
\newcommand\const{\mathrm{const}}
\newcommand\cS{\mathcal{S}}
\newcommand\cT{\mathcal{T}}
\newcommand\cU{\mathcal{U}}
\newcommand\defeq{\mathrel{\vcenter{\baselineskip0.5ex \lineskiplimit0pt \hbox{.}\hbox{.}}} =}
\newcommand\defequiv{\mathrel{\vcenter{\baselineskip0.7ex \lineskiplimit0pt \hbox{.}\hbox{.}}} \equiv}
\newcommand\defword[1]{\textbf{#1}}
\newcommand\enumbelow{\ }
\newcommand\ev{\mathrm{ev}}
\newcommand\Fam{\mathrm{Fam}}
\newcommand\FORALL{\text{ for all }}
\newcommand\Fun{\operatorname{Fun}}
\newcommand\FV{\mathrm{FV}}
\newcommand\Hom{\mathrm{Hom}}
\newcommand\hto{\hookrightarrow}
\newcommand\htox[2][]{\xhookrightarrow[#1]{#2}}
\newcommand\I{^{-1}}
\newcommand\id{\mathrm{id}}
\newcommand\im{\operatorname{im}}
\newcommand\inc{\mathrm{inc}}
\newcommand\loccit{\textit{loc.\@ cit.\@}}
\newcommand\Meas{\mathbf{Meas}}
\newcommand{\N}{\mathbb{N}}
\newcommand\nil{\mathsf{nil}}
\newcommand\norm[2][]{#1\lVert{#2}#1\rVert}
\newcommand\ol[1]{\overline{#1}}
\newcommand\op{\mathrm{op}}
\newcommand\opcit{\textit{op.\@ cit.\@}}
\newcommand\ot{\leftarrow}
\newcommand\pb[1][rd]{\ar[#1, phantom, "\lrcorner", pos=0]}
\newcommand\pBig[1]{\Big({#1}\Big)}
\newcommand\pbig[1]{\big({#1}\big)}
\renewcommand\phi{\varphi}
\newcommand\pow{\mathcal{P}}
\newcommand\presem[1]{[{#1}]}
\newcommand\R{\mathbb{R}}
\newcommand\reg{\mathrm{reg}}
\newcommand\rstr[2]{{\left.#1\right|_{#2}}}
\newcommand\SAn{\mathbf{SAn}}
\newcommand\SBor{\mathbf{SBor}}
\newcommand\sem[1]{[\![#1]\!]}
\newcommand\Set{\mathbf{Set}}
\newcommand\set[2][]{#1\{{#2}#1\}}
\newcommand\Sh{\mathrm{Sh}}
\newcommand\SOA{{2\mathrm{OA}}}
\newcommand\Sub{\operatorname{Sub}}
\newcommand\tildei{\tilde\imath}
\newcommand\tm{\mathbf{1}}
\newcommand\To{\Rightarrow}
\newcommand\TOA{{3\mathrm{OA}}}
\newcommand\toi{\xrightarrow{_\sim}}
\newcommand\toix[1]{\xrightarrow[{\raisebox{3.5pt}[0pt][0pt]{\ensuremath{\scriptstyle{\sim}}}}]{#1}}
\newcommand\tolon{{\,:\,}}
\newcommand\too{\twoheadrightarrow}
\newcommand\tot{\leftrightarrow}
\newcommand\toto{\rightrightarrows}
\newcommand\tox[2][]{\xrightarrow[#1]{#2}}
\newcommand\tto{\rightarrowtail}
\newcommand\two{\mathbf{2}}
\newcommand\uBo{\mathsf{uBo}}
\newcommand\und[1]{\underline{#1}}
\newcommand\unin{{\in}}
\newcommand\va{\vec{a}}
\newcommand\vW{\vec{W}}
\newcommand\vw{\vec{w}}
\newcommand\vX{\vec{X}}
\newcommand\vx{\vec{x}}
\newcommand\wh[1]{\widehat{#1}}
\newcommand\wt[1]{\widetilde{#1}}
\newcommand\xot[2][]{\xleftarrow[#1]{#2}}

\newcommand\concat{\string^}

\newcommand\uhr{{\upharpoonright}}
\newcommand\bbP{\mathbb{P}}
\usepackage{bbm}
\newcommand\bbone{\mathbbm{1}}
\newcommand\supp{\operatorname{supp}}
\newcommand\Cantor{\two^\N}

\ifcomments
\makeatletter
\@mparswitchfalse
\makeatother
\normalmarginpar
\paperwidth=\dimexpr \paperwidth + 6cm\relax
\marginparwidth=\dimexpr \marginparwidth + 5cm\relax

\newcommand{\comm}[1]{\todo[color=purple!50!white]{#1}}

\else
\newcommand{\comm}[1]{}
\fi

\DeclareFontFamily{U}{min}{}
\DeclareFontShape{U}{min}{m}{n}{<-> udmj30}{}
\newcommand\yon{\!\text{\usefont{U}{min}{m}{n}\symbol{'210}}\!}
\newcommand\ku{{\!\text{\usefont{U}{min}{m}{n}\symbol{'117}}\!}}

\theoremstyle{definition}
\newtheorem{defn}{Definition}[section]
\newtheorem{rmk}[defn]{Remark}

\theoremstyle{remark}

\theoremstyle{plain}
\newtheorem{propn}[defn]{Proposition}
\newtheorem{thm}[defn]{Theorem}
\newtheorem{cor}[defn]{Corollary}
\newtheorem*{cor*}{Corollary}
\newtheorem{lem}[defn]{Lemma}
\newtheorem*{claim}{Claim}

\let\oldproofname=\proofname
\renewcommand{\proofname}{\rm\bf{\oldproofname}}

\title{The Borel Topos}
\author{Joseph Helfer, Patrick Lutz, and Forte Shinko}
\date{}

\begin{document}
\begin{abstract}
We construct a topos for Borel mathematics: statements which are true in Borel mathemaics correspond in a precise way to statements which are true in the internal logic of this topos.
We investigate several other properties of this topos, including proving that it satisfies a Borel analogue of the Church-Turing thesis: the statement that every function on a standard Borel space is Borel.
We also show that an earlier construction of a similar ``Borel topos'' by van Osdol does not suit these purposes because of an error in that paper.
\end{abstract}

\maketitle
\tableofcontents

\comm{The right margin is extra wide as long as comments are on.}

\section{Introduction}
Descriptive set theory, or more precisely \emph{Borel mathematics}, is the study of the properties of mathematical objects built out of \emph{standard Borel spaces} (see \S\ref{subsec:std-borel})
using only Borel operations and relations.
For example, Borel graph combinatorics is the study of graphs whose set of vertices form a standard Borel space and whose set of edges is a Borel subset of that space (see \S\ref{subsec:borel-examples} for further examples).

One reason for interest in Borel mathematics is that on the one hand, Borel sets have many useful properties (for example, all Borel subsets of $\R$ are Lebesgue measurable) and on the other hand, many of the mathematical constructions which one ordinarily performs do not take one outside the world of Borel mathematics.

The apparent robustness of Borel constructions suggests the existence of a self-contained ``Borel mathematical universe'' in which \emph{all} mathematical constructions are in fact Borel constructions.
We are thus led to the following two related questions:
\begin{enumerate}[1.]
\item Can we find a model for set theory which satisfies precisely those statements which are true in Borel mathematics?
\item Can we (at least partially) characterize the theory of this model?
  I.e., can we find a list of axioms which exactly characterize the constructions which are allowed in Borel mathematics?
\end{enumerate}
The purpose of this paper is to give an answer to (1).
We also make some progress towards (2), although---despite this question being a large part of our motivation---we do not provide a definitive answer to it.

Our answer to question (1) is given by constructing an appropriate \emph{topos} (which is essentially the same thing as a model of higher-order arithmetic, see \cite{lambek-scott}) which can be thought of as the topos of Borel mathematics.
There are several reasons to use topoi as opposed to a model of a membership-based set theory (such as ZF) for this purpose.

The first is that for the purposes of Borel mathematics, the relevant aspects of sets is their role as carriers of abstract structures, and the mappings between them---rather than the membership relation and the cumulative hierarchy.
In other words, what is relevant is \emph{abstract set theory} in the sense of \cite{makkai-abstract-sets}, which is exactly that aspect of set theory which is captured by topos theory.
Secondly, there is already a well-established construction in topos theory which is precisely suited to our purpose of beginning with a basic category of interest (here, the standard Borel spaces), and expanding it into a model of set theory (i.e., a topos).
This is the construction of the \emph{category of sheaves} on a \emph{Grothendieck site}, which we review in \S\ref{subsection:topologies}.
Previous examples of such constructions in a similar spirit are Johnstone's \emph{topological topos} \cite{johnstone-topological-topos} and Dubuc's \emph{smooth topos} \cite{moerdijk-reyes}---and in fact, there is even a previously existing \emph{Borel topos} \cite{van-osdol-borel-topos}, though as we explain in \S\ref{subsec:van-osdol-topos}, this topos does not suit our purposes.
(Also related is the \emph{effective topos} of \cite{hyland-effective-topos}, though it is constructed using a different, but related, method.)

A third reason, and one that is at first sight not obviously relevant, is the fact that a topos is in general a model only of \emph{intuitionistic} set theory.
Indeed, the Borel topos we construct does not obey classical logic (i.e., it is not a \emph{Boolean topos}), and for reasons we will now explain, this turns out to work in our favour.
That there is a connection between Borel mathematics and constructivism is also not a new idea; cf.\ \cite{shani-borel-reducibility-symmetric-models,bernshteyn-descriptive-distributed-survey,bernshteyn-distributed-lovasz}.
\comm{
  More details about this:
  Shani has found a close connection between certain statements involving Borel equivalence relations and implications between weak forms of the Axiom of Choice~\cite{shani-borel-reducibility-symmetric-models} and Bernstein has discovered a close connection between the existence of Borel solutions to certain problems in graph combinatorics and the existence of a fast distributed algorithm for solving the corresponding problem on finite graphs~\cite{bernshteyn-descriptive-distributed-survey,bernshteyn-distributed-lovasz}. Our construction of the topos of Borel mathematics can be seen as another way to make the connection between Borel mathematics and constructive reasoning precise.
}

The relevance of intuitionistic logic to the Borel topos is related to the infeasibility of creating a universe which truly only consists of Borel sets.
Rather, one must also admit in the \emph{analytic sets}: the images of Borel sets under Borel maps.
This is basically unavoidable, since in any topos, we must be able to form the image of an arbitrary map.
Having admitted the analytic sets into our universe, it is then natural to ask for an intrinsic description of those subsets which are in fact Borel, and this turns out to have a very satisfying answer: the subsets \( P \) of a set \( X \) which are Borel are precisely those which are \emph{decidable} in the sense that \( x \in P \vee x \notin P \) for all \( x \in X \).

Let us say this somewhat more precisely: each standard Borel set \( X \) gives rise to a ``set'' (i.e., an object) \( X \) in our topos, and the subobjects of \( X \) correspond to the analytic subsets of \( X \) (more precisely, to the \emph{\( \sigma \)-ideals} of analytic sets, see \S\ref{subsec:subsheaves-in-e}), and those subobjects which---externally---correspond to Borel subsets of \( X \) are---internally---precisely the decidable subobjects of \( X \).
But in fact, we can say more: it also makes sense to ask \emph{internally} whether it is true that ``a subset (say, of \( \two^\N \)) is decidable if and only if it is Borel'', and we prove in \S\ref{subsec:borel-myopia} that this is so.
Moreover, as we explain there, this is intuitionistically equivalent to the statement ``every function is Borel'', hence the latter holds in our topos as well.
This is a Borel analogue to a strong version of the Church-Turing thesis, stating that ``every function is computable'', which can likewise be made to hold literally inside of an appropriate topos.

\subsection{Main results}
Our main theorem expresses that our topos \( \cE \) is indeed ``the topos of Borel mathematics'' in the sense that a statement \( \phi \) holds internally in \( \cE \) if and only if the ``Borel version'' of \( \cE \) is true.
To make this precise, we need to say what kind of statements \( \phi \) we are considering, and what we mean by the Borel version of such a statement.
We take \( \phi \) to be a statement of \emph{third order arithmetic} in the sense described in \S\ref{subsec:internal-external}: the first-order language with variables of sorts \( \N \), \( \two^\N \), and \( \two^{\two^\N} \).
The ``Borel interpretation'' of such a statement is just the usual, set-theoretic interpretation, except that the variables of sort \( \two^{\two^\N} \) are interpreted as quantifying over \emph{Borel} subsets of \( \two^\N \).
There is an additional implicit informal claim here that needs to be justified, namely that all claims of ``Borel mathematics'' \emph{are} the Borel interpretation of some such third-order statement.
We provide evidence for this in \S\ref{subsec:borel-examples} by showing that a variety of statements in the literature are indeed of this form.

There is an additional twist involved in the statement of this theorem, which is that the Borel interpretation is sound with respect to classical logic whereas the interpretation of sentences in our topos \( \cE \) is, a priori, only sound for intuitionistic logic.
And in fact, all that is true in general is that for each sentence \( \phi \), there is a \emph{classically equivalent} statement \( \tilde \phi \) such that the Borel interpretation of \( \phi \) holds if and only if \( \tilde \phi \) holds in \( \cE \).
The particular \( \tilde \phi \) we take is the \emph{Kuroda translation} of \( \phi \), which is one of several versions of the ``double negation interpretations'' of classical into intuitionistic logic.
We thus arrive at
{\renewcommand\thedefn{A}
\begin{thm}[see Theorem~\ref{thm:internal-external}]
  For any sentence \( \phi \) in third-order arithmetic, the Borel interpretation of \( \phi \) holds if and only if its Kuroda translation \( \tilde \phi \) holds in the Borel topos \( \cE \).
\end{thm}
}
Since topoi soundly interpret all of higher-order arithmetic (say, in the sense of \cite{lambek-scott}), we thus obtain the following criterion for the truth of statements in Borel mathematics:
\begin{cor*}
  For any sentence \( \phi \) of third-order arithmetic, if the Kuroda translation \( \tilde \phi \) is provable in intuitionistic higher-order arithmetic, then the Borel interpretation of \( \phi \) is true.
\end{cor*}
In fact, more than this is true, since it is not only pure higher-order arithmetic which is valid in our topos \( \cE \).
For example, in proving \( \tilde \phi \), we may additionally assume the Kuroda translations of any statements whose Borel interpretations have already been established.

The second main result, which we call \emph{Borel myopia} is the one already mentioned above:
{\renewcommand\thedefn{B}
\begin{thm}[see \S\ref{subsec:borel-myopia} for a precise statement]
  The following statement holds internally in the Borel topos \( \cE \): a subset of \( \two^\N \) is decidable if and only if it is Borel.
\end{thm}
}

\subsection{Further results}
In addition to the above main results, we prove a number of other basic topos-theoretic properties of our Borel topos \( \cE \):
\begin{itemize}
\item In \S\ref{subsec:more-topos-props}, we prove that \( \cE \) satisfies Markov's principle, and that it is two-valued but (as already noted) not Boolean, and hence is not well-pointed and does not satisfy the axiom of choice.
\item In \S\ref{subsec:subsheaves-in-e}, we show that the subobject lattice of \( \two^\N \) in \( \cE \) is isomorphic to the lattice of \( \sigma \)-ideals of analytic subsets of \( \two^\N \).
\item In \S\ref{subsec:double-neg-subtopos}, we show that the double-negation subtopos of \( \cE \) is equivalent to \( \Set \).
\item In \S\ref{subsec:can-top}, we show that, assuming the continuum hypothesis, \( \cE \) is equivalent to the category of sheaves with respect to the canonical topology on the category \( \SBor \) of standard Borel spaces.
In Appendix~\ref{sec:necessity-of-ch}, we present an argument of Edward Hou showing that this is not provable in ZFC (without CH).
\end{itemize}

Moreover, in \S\ref{subsec:analytic}, we show that \( \cE \) can be described as sheaves on the category of \emph{analytic} spaces (i.e., analytic subsets of complete Borel spaces), rather than complete Borel spaces.
We also prove various other properties of this category, and in particular show that it is the \emph{completion}, in a certain sense---and in fact in several senses---of the category of standard Borel spaces.

\subsection{Statement on lack of AI usage.}
No AI tools were used in the creation of paper.

\section{Preliminaries}
We will assume familiarity with the basic concepts of category theory and of descriptive set theory, for which we refer to \cite{maclane-categories} and \cite{kechris-classical}, respectively.
For additional background on topoi beyond what we review here, we refer to \cite{johnstone-topos-theory,mac-lane-moerdijk,johnstone-elephant-vols}.

\subsection{Grothendieck topologies}\label{subsection:topologies}
We recall some basic notions and fix our notation regarding Grothendieck topologies.

Let $\cC$ be a small category with finite limits.
We recall that there are two ways of defining a Grothendieck topology on $\cC$: in terms of \emph{sieves}, or in terms of \emph{covering families} (see \cite[\S1.1]{makkai-reyes-focl}).
We will use the latter, though we will return to sieves in \S\S\ref{subsec:subsheaves}~and~\ref{subsec:more-topos-props}.

Given an object $X \in \cC$, a \defword{family of morphisms over} $X$ is a set $\bA \subset \Ar(\cC)$ of morphisms with codomain $X$.
We write $\Fam(X)$ for the set of families over $X$.
For \( \bA \in \Fam(X) \) and \( a \in \bA \), we typically write \( U_a \) for the domain of \( a \).
A \defword{(Grothendieck) topology} $J$ on $X$ consists of a set $J(X) \subset \Fam(X)$ of families of morphisms over $X$, called the \defword{covering families}, or \defword{coverings} in $J$, for each $X \in \cC$, subject to the following axioms:

(T1) The singleton family of morphisms $\set{\id_X}$ is in $J(X)$ for each $X \in \cC$

(T2) Given a covering family $\bA \in J(X)$, a morphism $f \colon Y \to X$ in $\cC$, and a choice of pullback
\[
  \begin{tikzcd}
    f^*U_a \ar[r, "\pi_f"] \ar[d, "f^*a"'] \pb & U_a \ar[d, "a"] \\
    Y \ar[r, "f"] & X
  \end{tikzcd}
\]
for each $a \in \bA$, the resulting family $f^*\bA \defeq \set{f^*U_a \tox{f^*a} Y}_{a \in \bA}$ of morphisms over $Y$ is in $J(Y)$.

(T3) Given a covering family $\bA \in J(X)$ and further covering families $\bB_a \in J(U_a)$ for each $a \in \bA$, the resulting family of morphisms $\set{a \circ b \colon U_b \to X \mid a \in \bA, b \in \bB_a }$ is also in $J(X)$.

(T4) Given $\bA \in J(X)$ and $\bB \in \Fam(X)$, if each $a \in \bA$ factors through some $b \in \bB$ (i.e., there exists $g \colon U_a \to U_b$ with $b \circ g = a$), then $\bB \in J(X)$.

A category with finite limits together with a Grothendieck topology is called a \defword{site}.
We may conflate a site with its underlying category.

Given a presheaf $F$ on $\cC$ (i.e., a functor $F \colon \cC^\op \to \Set$) and a family of morphisms $\bA \in \Fam(X)$ in $\cC$, we say that $F$ is a \defword{sheaf with respect to $\bA$}, or an \defword{$\bA$-sheaf}, and that $\bA$ is an \defword{\boldmath $F$-good family of morphisms}, if given pullbacks
\[
  \begin{tikzcd}
    U_a \times_X U_{a'} \ar[r, "\pi"] \ar[d, "\pi'"'] \pb & U_a \ar[d, "a"] \\
    U_{a'} \ar[r, "a'"] & X
  \end{tikzcd}
\]
for all $a,a' \in \bA$, and elements $s_a \in F(U_a)$ for each $a \in \bA$ satisfying $\pi^*s_a = (\pi')^*(s_{a'}) \in F(U_a \times_X U_{a'})$ for all $a,a' \in \bA$, there exists a unique $s \in F(X)$ with $s_a = a^*s$ for all $a \in \bA$.
Given a topology $J$ on $\cC$, we say that a presheaf $F$ is a \defword{sheaf with respect to $J$} (or \defword{for $J$}, or a \defword{sheaf on $\cC$}) if it is one with respect to $\bA$ for each $X \in \cC$ and $\bA \in J(X)$.

We can just as well speak about $F$ being a sheaf with respect to \emph{any} collection of families of morphisms in $\cC$, not just a Grothendieck topology.
The main significance of the notion of Grothendieck topology is that if two topologies $J$ and $J'$ on $\cC$ have the same sheaves (i.e., $F \colon \cC^\op \to \Set$ is a sheaf for $J$ if and only if it is a sheaf for $J'$), then they are equal.
More precisely, given a topology $J$, the covering families $J(X) \subset \Fam(X)$ are exactly those $\bA \in \Fam(X)$ such that $f^*\bA$ is $F$-good for every $J$-sheaf $F$ and every $f \colon Y \to X$ in $\cC$.

Given any set $J$ of families of morphisms in $\cC$, there is a least topology containing $J$, called the topology \defword{generated by} $J$, and which we denote by $\br{J}$.
If $J$ satisfies the axiom (T2), then the sheaves for $\br{J}$ are the same as those for $J$.
It follows that two sets $J,J'$ of families of morphisms, both satisfying (T2), generate the same topology if and only if they have the same sheaves.

If $J$ moreover satisfies (T1)-(T3), it is called a \defword{pretopology}.
The significance of this notion is that in this case,
\comm{
  (and in fact, even if (T1) is replaced by ``the set $J(X)$ is inhabited for each $X \in \cC$'')%
}
we can say precisely what the topology generated by $J$ is: $\br{J}(X)$ consists precisely of those $\bB \in \Fam(X)$ such that some $\bA \in J(X)$ \defword{factors} through $\bB$, by which we mean that each $a \in \bA$ factors through some $b \in \bB$, as in (T4).
We call a category with finite limits together with a pretopology a \defword{presite}.

For a category $\cC$, we write $\wh{\cC} = \Fun(\cC^\op,\Set)$ for the category of presheaves on $\cC$, and given a topology $J$ on $\cC$, we write $\Sh(\cC,J) \subset \wh{\cC}$ (or just $\Sh(\cC)$) for the category of $J$-sheaves.
A \defword{Grothendieck topos} (or just \defword{topos}) is a category which is equivalent to $\Sh(\cC,J)$ for some site $(\cC,J)$.

We write $\yon \colon \cC \to \wh{\cC}$ for the Yoneda embedding $\yon(X) = \cC(-,X)$, and we also write $\wh{X}$ for $\yon(X)$; we recall that
the Yoneda embedding preserves all small limits.
A presheaf is \defword{representable} if it is isomorphic to one of the form \( \wh{X} \).
The topology \( J \) is called \defword{subcanonical} if all representable presheaves are sheaves, i.e., \( \yon \) factors through \( \Sh(\cC) \).
The topologies we consider will all be subcanonical.

The full subcategory \( \Sh(\cC) \subset \wh{\cC} \) is \emph{reflective}, meaning that the inclusion \( i \colon \Sh(\cC) \hto \wh{\cC} \) has a left adjoint \( \wh{\cC} \to \Sh(\cC) \), called \defword{sheafification}.
In particular, \( i \) preserves and reflects all small limits, and hence also monomorphisms, and hence \( \yon \) also preserves monomorphisms.

\begin{lem} \label{lem:empty-object}
  If \( (\cC, J) \) is a small subcanonical presite, then an object \( U \in \cC \) is covered by the empty family if and only if \( \wh U \in \Sh(\cC, J) \) is initial.
\end{lem}
\begin{proof}
  The forward implication follows immediately from the definitions and Yoneda's lemma.
  For the reverse implication, suppose \( \wh U \) is initial, and let us write \( \emptyset_U \in \Fam(U) \) for the empty cover of \( U \).
  Then, for any morphism \( f \colon U' \to U \) in \( \cC \), the morphism \( \wh f \colon \wh{U'} \to \wh U \) is an isomorphism, and \( \wh{U'} \) is likewise initial (this is a general property of topoi, and more generally, of cartesian closed categories).
  It follows that \( \emptyset_{U'} = f^* \emptyset_U \) is \( F \)-good for any \( F \in \Sh(\cC, J) \), so we conclude that \( \emptyset_U \in \br{J}(U) \).
  But this means that there is some cover \( \bA \in J(U) \) which factors through \( \emptyset_U \).
  Hence, we must have \( \emptyset_U \in J(U) \), since the only cover factoring through \( \emptyset_U \) is \( \emptyset_U \) itself.
\end{proof}

\subsection{Subsheaves}\label{subsec:subsheaves}
In any presheaf category \( \wh{\cC} \), the subobjects of a given presheaf \( F \in \wh{\cC} \) can be described explicitly.
Namely, a \defword{subfunctor} of \( F \) is a presheaf \( R \in \wh{\cC} \) such that \( R X \subset F X \) for \( X \in \cC \), and for which \( R f \colon R Y \to R X \) is the restriction of \( F f \colon F Y \to F X \) for each \( f \colon X \to Y \) in \( \cC \); in particular, \( (F f) (R Y) \subset R X \).
Conversely, an arbitrary prescription of subsets \( R X \subset F X \) for \( X \in \cC \) satisfying \( (F f) (R Y) \subset R X \) determines a subfunctor.
The inclusion maps \( R X \hto F X \) then give a natural transformation \( R \hto F \), which is a monomorphism in \( \wh{\cC} \).

The fact is now that, every monomorphism \( i \colon G \to F \) in \( \wh{\cC} \) is isomorphic to a unique subfunctor \( R \) of \( F \), namely with \( R(X) = i_X(G X) \) for \( X \in \cC \).

A subfunctor \( R \) of a representable \( \wh{X} \) is called a \defword{sieve} on \( X \); concretely, it is the specification of a subset \( R Y \subset \cC(Y,X) \) for each \( Y \in \cC \) such that, given \( Z \tox{g} Y \tox{f} X \), if \( f \in R Y \), then \( f \circ g \in R Z \).

Now fix a topology \( J \) on \( \cC \), and a sheaf \( F \in \Sh(\cC) \).
Because the inclusion \( \Sh(\cC) \hto \wh{\cC} \) preserves and reflects monomorphisms, the monomorphisms \( G \to F \) into \( F \) in \( \Sh(\cC) \) are precisely the monomorphisms in \( \wh{\cC} \) for which \( G \) is a sheaf, and hence are in bijection with those subfunctors \( R \subset F \) for which \( R \) is a sheaf, which are called \defword{subsheaves}.

The subsheaves of a representable \( \wh{X} \) are those sieves \( R \subset \wh{X} \) such that, given \( f \colon Y \to X \) and a covering \( \bA \in J(Y) \), if \( f \circ a \colon U_a \to X \) is in \( R A \) for each \( a \in \bA \), then \( f \in R Y \).

The subsheaves of any sheaf \( F \), ordered by inclusion, form a Heyting algebra \( \Sub(F) \), whose operations are determined by the KJ semantics discussed in \S\ref{subsec:kj-semantics} below.
However, we note that finite conjunctions in \( \Sub(F) \) are simply given pointwise: \( (R \wedge S) X = R X \cap S X \).

\subsection{Elementary topoi} \label{subsec:elementary-topoi}
Though we are mainly concerned with Grothendieck topoi in this paper, we will also have occasion to consider the more general \emph{elementary topoi}.

We recall that a \defword{power object} of an object \( X \) in a category \( \cD \) is an object \( \pow(X) \) equipped with a binary relation \( \in_X \tto X \times \pow(X) \) with the universal property that any other binary relation \( R \tto X \times U \) is pulled back along \( \id \times f \colon X \times U \to X \times \pow X \) for a unique morphism \( f \colon U \to \pow X \) (which we say \emph{classifies} \( R \)).

An \defword{elementary topos} is a category with finite limits such that every object admits a power object.
Every Grothendieck topos is an elementary topos, and many of the most important properties of Grothendieck topoi are shared by general elementary topoi.
For one, elementary topoi have all finite colimits, and exponential objects.
But most importantly, any elementary topos \( \cD \) has the following \emph{logical} structure: each subobject lattice \( \Sub(X) \) for \( X \in \cD \) is a Heyting algebra, and for each morphism \( f \colon X \to Y \) in \( \cD \), the pullback morphism \( f^* \colon \Sub(Y) \to \Sub(X) \) is a homomorphism of Heyting algebras and has left and right adjoints \( \exists_f,\forall_f \colon \Sub(X) \to \Sub(Y) \).

This gives rise to the \emph{internal logic} of \( \cD \).
Consider the many-sorted first-order language \( \cL_\cD \) whose sorts are the objects of \( \cD \), with a function symbol (with evident input and output sorts) for each morphism \( f \colon \prod_{i=1}^n X_i \to X \) out of a product in \( \cD \), and a relation symbol (again with evident arity) for each monomorphism \( R \tto \prod_{i=1}^n X_i \) into a product in \( \cD \).
We write \( x \tolon X \) to indicate that \( x \) is a variable of sort \( X \), and given sequences \( \vec{x} = (x_1,\ldots,x_n) \) and \( \vX = (X_1,\ldots,X_n) \), we write \( \vec{x} \tolon \vX \) to indicate that \( x_i \tolon X_i \) for each \( i \).

Using the aforementioned logical structure, we may then, as in \cite[Chapter~2]{makkai-reyes-focl}, define, for each term \( t \) in \( \cL_\cD \) with free variables \( \FV(t) \subset \vx \tolon \vX \), an arrow \( \presem{t}_{\vx} \colon \prod_{i=1}^n X_i \to X \), and for each formula \( \phi \) of \( \cL_\cD \) with \( \FV(\phi) \subset \vx \tolon \vX \) a subobject \( \presem{\phi}_{\vx} \in \Sub(\prod_{i=1}^n X_i) \).
In particular, the propositional connectives are interpreted using the Heyting algebra structure on \( \Sub(\prod_{i=1}^n X_i) \).
If \( \phi \) is a closed formula, \( \presem{\phi}_{\emptyset} \) is thus a subsheaf of the terminal sheaf \( \tm \), and we say \( \cD \vDash \phi \) if \( \presem{\phi}_{\emptyset} = \top_{\Sub(\tm)} \).

We now make some remarks about the soundness of this interpretation and various extensions of it:
\begin{enumerate}[1.]
\item The fundamental fact about this interpretation is that it is sound for intuitionistic logic, in the sense that \( \cD \vDash \phi \) for every intuitionistic validity \( \phi \).
\item Moreover, it is sound for intuitionistic \emph{higher-order} logic, in the sense that, if we interpret second-and-higher-order variables using the power objects in \( \cD \), then the axioms of higher-order logic (extensionality and comprehension) are satisfied as well.\footnote{
    Technically, this extended interpretation depends on the \emph{choice} of power objects for each object in \( \cD \), as these are only uniquely defined up to isomorphism.
    However, the \emph{truth values} in the resulting interpretation are of course independent of these choices.
    Similar remarks apply to the further extended interpretations described below.
  }
\item Recall that a \defword{natural numbers object} (NNO) in a category is an object \( N \) equipped with a universal diagram of the shape \( \tm \tox{0} N \tox{s} N \) with \( \tm \) terminal: given any other diagram \( \tm \tox{a} A \tox{f} A \), there is a unique morphism \( F \colon N \to A \) with \( F 0 = a \) and \( F s = f F \).

  If the elementary topos \( \cD \) contains an NNO, then by interpreting the natural numbers using the NNO and interpreting higher-order sorts using power objects as above, we obtain a sound interpretation of intuitionistic \emph{higher-order arithmetic} (HOA): the theory with sorts \( \pow^n(\N) \) for \( n \ge 1 \), function symbols \( 0 \tolon \N \) and \( s \colon \N \to \N \), and relation symbols \( {\in} \tolon \pow^n(\N) \times \pow^{n+1}(\N) \) for each \( n \) and with axioms extensionality and comprehension for each \( n \) and the (second-order) Peano axioms.
\item Finally, this interpretation can be extended to a sound interpretation of any \emph{definitional extension} of HOA, the meaning of which we now explain.

  By a \defword{formula-in-context}, we mean a pair \( (\phi, \vx) \) consisting of a formula \( \phi \) and a sequence \( \vx \) of variables including the free variables of \( \phi \).
  By an \defword{immediate definitional extension} of a theory \( (\cL, \cT) \), we mean an extended theory \( (\cL', \cT') \) obtained from three sets \( (\Phi_{\mathrm{sort}}, \Phi_{\mathrm{fun}}, \Phi_{\mathrm{rel}})  \) of formulas-in-context over \( \cL \) as follows:
  \begin{itemize}
  \item \( \cL' \) has an additional sort \( A_\phi \) for each \( (\phi, \vx \tolon \vX) \in \Phi_{\mathrm{sort}} \) and function symbols \( \pi_i \colon A_{\phi} \to X_i \), and \( \cT' \) has an additional axiom stating that the \( \pi_i \) are jointly injective and identify \( A_{\phi} \) with the subset of \( \prod_{i} X_i \) satisfying \( \phi \).
  \item For each \( (\phi, \vx \colon \vX) \in \Phi_{\mathrm{fun}} \), \( \vx \) is required to be of the form \( \vw y \colon \vW Y \) (i.e., it must have at least one variable), and \( \phi \) is required to be, provably in \( \cT \), a function relation (with output \( y \)); \( \cL' \) then has an additional function symbol \( f_\phi \colon \vW \to Y \), and \( \cT' \) has the additional axiom \( f_\phi(\vw) = y \tot \phi \)
  \item \( \cL' \) has an additional relation \( R_\phi \tolon \vX \) for each \( (\phi, \vx \tolon \vX) \in \Phi_{\mathrm{rel}} \), and \( \cT' \) has the additional axiom \( R_\phi(\vx) \tot \phi \).
  \end{itemize}
  By a \defword{definitional extension} of \( (\cL, \cT) \), we mean an extension obtained from a finite sequence of immediate definitional extensions.

  The interpretation of HOA in an elementary topos with NNO extends in an evident manner to any definitional extension of HOA.
  The fact that this interpretation remains sound follows from the more basic fact that a definitional extension \( (\cL', \cT')  \) of any theory \( (\cL, \cT) \) is \emph{logically conservative}: any sentence in \( \cL \) provable in \( \cT' \) is already provable in \( \cT \).
  But more than this, a definitional extension is \emph{expressively conservative} in that any sentence of \( \cL' \) is, provably in \( \cT' \), equivalent to some sentence in \( \cL \).
  Thus, a definitional extension allows us neither to prove, nor to express, more theorems, but only to express things more concisely.

  We will exhibit a particular definitional extension of HOA in \S\ref{subsec:borel-myopia} which is suited to our purposes, and in \S\ref{subsec:borel-examples}, we will consider a variant of the notion of definitional extension.
\end{enumerate}

\subsection{Kripke-Joyal semantics}\label{subsec:kj-semantics}
In case our elementary topos \( \cD \) is the category of sheaves \( \cD = \Sh(\cC) \) on a presite \( (\cC,J) \) (which for simplicity we assume to be subcanonical), the internal logic of \( \cD \) takes on an especially clear and explicit form, known as the \emph{KJ semantics} or \emph{sheaf semantics}.\footnote{%
  There is also a version of the KJ semantics which makes sense in a general elementary topos.
}

We define the following ``forcing relation'': given an object \( U \in \cC \), morphisms \( a_i \colon \wh U \to X_i \), and a formula \( \phi \) in \( \cL_{\cD} \) with \( \FV(\phi) \subset \vx \tolon \vX \), we write \( U \Vdash_{\vx \mapsto \va} \phi \) (or just \( U \Vdash \phi \) if \( \vx \) and \( \va \) are clear from context) if the morphism \( \br{a_i}_{i=1}^n \colon \wh U \to \prod_{i=1}^n X_i \) factors through \( [\phi]_{\vx} \hto \prod_{i=1}^n X_i \).

We now list some crucial properties satisfied by the forcing relation (see \cite[\S\S VI.6-7]{mac-lane-moerdijk} for proofs).
For a term \( t \) of sort \( A \) with \( \FV(t) \subset \vx \tolon \vX \) and morphisms \( a_i \colon \wh U \to X_i \), let us write \( \sem{t}_{\vx,\va} \) (or just \( \sem{t} \)) for the morphism \( [t]_{\vx} \circ \br{a_i}_{i=1}^n \colon \wh{U} \to A \).
We then have:
\begin{propn}\label{propn:kripke-joyal}\enumbelow
  \begin{enumerate}[(i)]
  \item With \( U,\vx,\va \) as above, if \( U \Vdash \phi \), then \( U' \Vdash \phi \) for any morphism \( u \colon U' \to U \); note that by \( U' \Vdash \phi \), we mean \( U' \Vdash_{\vx \mapsto \va'} \phi \) where \( \va' = \br{a_1 \circ \yon(u),\ldots,a_n \circ \yon(u)} \).
  \item \( U \Vdash \phi \) whenever (the universal closure of) \( \phi \) is intuitionistically valid.
  \item If \( \bA \in J(U) \) is a covering and \( U_a \Vdash \phi \) for each \( a \in \bA \), then \( U \Vdash \phi \).
  \item If \( \phi \) is closed, then \( \cD \vDash \phi \) is equivalent to \( \forall\, U \in \cC.\ U \Vdash_{\emptyset,\emptyset} \phi \), and if \( \cC \) has a terminal object \( \tm \), then also to \( \tm \Vdash \phi \).
  \item\label{item:forcing-conditions} The relation \( U \Vdash_{\vx \mapsto \va} \phi \) is given explicitly by the following recursive clauses:
  \end{enumerate}
  {
    \def\arraystretch{1.1}
    \begin{center}
      \begin{tabular}{lp{11.5cm}}
        \( U \Vdash \top \) & always. \\
        \( U \Vdash \bot \) & if \( \wh U \) is an initial object. \\
        \( U \Vdash s = t \) & if \( \sem{s} = \sem{t} \). \\
        \( U \Vdash R t_1 \ldots t_m \) & (\( R \subset \prod_{i=1}^n X_i \) a relation) if \( \br{\sem{t_i}}_{i=1}^n \colon \wh U \to \prod_{i}X_i \) factors through \( R \). \\
        \( U \Vdash \phi_1 \wedge \phi_2 \) & if \( U \Vdash \phi_1 \) and \( U \Vdash \phi_2 \). \\
        \( U \Vdash \phi_1 \vee \phi_2 \) & if there is \( \bA \in J(U) \) such that for \( a \in \bA \), either \( U_a \Vdash \phi_1 \) or \( U_a \Vdash \phi_2 \). \\
        \( U \Vdash \phi_1\to\phi_2 \) & if \( U' \Vdash \phi_1 \) implies \( U' \Vdash \phi_2 \) for all \( U'\to U \) in \( \cC \). \\
        \( U \Vdash \forall x_{n+1}\,\phi \) & if \( U' \Vdash \phi \) for all \( u\colon U'\to U \) and all \( a_{n+1} \colon \wh{U'} \to X_{n+1} \), where \( x_{n+1} \tolon X_{n+1} \). \\
        \( U \Vdash \exists x_{n+1}\,\phi \) & if there is \( \bB \in J(U) \) and morphisms \( a_b \colon \wh{U_b} \to X_{n+1} \) for \( b \in \bB \) such that \( U_b \Vdash \phi \) for all \( b \in \bB \), where \( x_{n+1} \tolon X_{n+1} \).
      \end{tabular}
    \end{center}
  }
\end{propn}
\begin{rmk}
  The clauses for \( \bot \), \( \vee \), and \( \exists \) are often somewhat simpler in practice---and in particular in the topoi \( \cE \) and \( \cE^\inc \) introduced below in \S\ref{sec:borel-topos}.
  Namely, we have:
  \begin{itemize}
  \item \( U \Vdash \bot \) iff \( U \) is initial.\\
    This holds as long as an object in \( \cC \) is covered by the empty family iff it is initial (by Lemma~\ref{lem:empty-object}).
  \item \( U \Vdash \phi \vee \psi \) iff there is a covering \( B_1,B_2 \to U \) such that \( B_i \Vdash \phi_i \) for \( i = 1,2 \).\\
    This holds assuming that for any \( U \in \cC \) and \( \bA_1,\bA_2 \in \Fam(X) \) such that \( \bA_1 \cup \bA_2 \) is a cover, there are morphisms \( A_1,A_2 \to U \) such that each \( U_a \tox{a} U \) in \( \bA_i \) factors through \( A_i \); for our topoi \( \cE \) and \( \cE^\inc \), this condition follows from the existence of finite coproducts in the site \( \SBor \).
  \item \( U \Vdash \exists x_{n+1}\,\phi \) iff there is a cover \( U' \to U \) and \( a \colon U' \to X_{n+1} \) such that \( U' \Vdash \phi \).\\
    This holds for \( \cE \)---again using the existence of countable coproducts in \( \SBor \)---but does \emph{not} hold for \( \cE^\inc \).
    \comm{
      Proof (for \( \cE) \)):
      This follows since, if \( \bB \in J_{\SBor}(U) \) is any covering, then so is the singleton family \( U' = \coprod_{b \in \bB} U_b \to U \), and the coproduct \( U' = \coprod_{b \in \bB} U_b \) is preserved by \( \yon \colon J_{\SBor} \to \cE \).
      Then, given \( a_b \colon U_b \to X_{n+1} \) with \( U_b \Vdash \phi \) for each \( b \in \bB \), we obtain \( a \colon U' \to X_{n+1} \), and since \( U_b \Vdash \phi \) for each \( b \) in the cover \( \set{U_b \tox{\mathrm{in}_b} U'}_{b \in \bB} \), it follows that \( U' \Vdash \phi \).
    }
  \end{itemize}
\end{rmk}

\subsection{Standard Borel spaces} \label{subsec:std-borel}
A \defword{Polish space} $X$ is a separable complete metrizable topological space.
A \defword{standard Borel space} is the underlying measurable space of a Polish space \( X \), i.e., it is \( X \) equipped with its \( \sigma \)-algebra of \defword{Borel subsets}: the least \( \sigma \)-algebra on \( X \) containing the open subsets.
We recall also that an \defword{analytic subset} of a standard Borel space \( X \) is any subset which is the image of some Borel (i.e., measurable) map \( Y \to X \) from a standard Borel space \( Y \).

Let $\SBor$ be the category of standard Borel spaces with Borel maps.
It is a full subcategory of the category \( \Meas \) of measurable spaces, and it is essentially small, each object being isomorphic to a finite set or the natural numbers $\N$ with the discrete topology, or to the real numbers $\R$ with the Euclidean topology (or equivalently, to $\two^\N$ or $\N^\N$ with the product topology) \cite[Theorem~15.6]{kechris-classical}.

The category \( \SBor \) has countable%
\footnote{%
  By ``countable'' we always mean ``finite or countably infinite''.%
}
 limits and coproducts \cite{chen-universal-characterization-of-standard-borel},
\comm{
  The diagram $f, g \colon X \rightrightarrows Y$ has a coequalizer
  iff quotient of $Y$ by the equivalence relation
  generated by $(f(x), g(x))_{x \in X}$ is standard Borel.
}
which are characterized by:
\begin{propn}\label{propn:sbor-lims-colims}
  The underlying-set functor \( U \colon \SBor \to \Set \) preserves and reflects countable limits and countable coproducts.
\end{propn}

We recall that this means that for a countable diagram \( D \colon J \to \SBor \), a given cone \( \ol D \colon J^\con \to \SBor \) over \( D \) is a limit cone if and only if \( U \circ \ol D \colon J^\con \to \Set \) is a limit cone in \( \Set \), and likewise with coproducts.

\begin{proof}
  \( U \) is clearly \emph{conservative}, meaning \( f \colon X \to Y \) in \( \SBor \) is an isomorphism if and only if \( U f \) is.
  Since \( \SBor \) has all countable limits and coproducts, it follows on general grounds that if \( U \) preserves countable limits and colimits, then it also reflects them.
  It thus remains to show the former.

  That \( U \) preserves limits follows from its being isomorphic to the hom-functor \( \SBor(\tm,-) \) from the terminal object \( \tm \in \SBor \).
  That \( U \) preserves countable coproducts is easy to see directly, since given a countable family \( (X_i)_{i \in I} \) of Polish spaces, a coproduct \( \coprod_{i \in I} X_i \) in \( \SBor \) can be constructed by equipping the disjoint union of the sets \( X_i \) with the usual topology.
\end{proof}

Next, we recall that given objects \( X,Y \in \cC \) in a category \( \cC \) with binary products, an \emph{exponential object} of \( Y \) by \( X \) is an object \( Y^X \) equipped with a morphism \( e \colon Y^X \times X \to Y \) such that for each \( A \in \cC \), the composite \( \cC(A,Y^X) \tox{\times X} \cC(A \times X, Y^X \times X) \tox{\ev \circ} \cC(A \times X, Y) \) is a bijection.

Given \( X,Y \in \SBor \) with \( X \) countable (hence discrete), we have the product \( Y^X = \prod_{x \in X} Y \in \SBor \) and the evaluation map \( \ev \colon Y^X \times X \to Y \), which is Borel, since for \( U \subset Y \) Borel, \( \ev\I(U) = \bigcup_{x \in X} \pi_x\I(U) \times \set{x} \), where \( \pi_x \colon \Pi_{x \in X} Y \to Y \) is the projection.

\begin{propn}\label{propn:sbor-exponentials}
  \( \ev \) exhibits \( Y^X \) as an exponential object of \( Y \) by \( X \) in \( \SBor \).
\end{propn}

\begin{proof}
  After unfolding the definitions, this reduces to the claim that, for each \( A \in \SBor \) and Borel map \( f \colon A \times X \to Y \), the map \( g \colon A \to Y^X \) with \( g(a)(x) = g(a,x) \) is Borel.
  However, since each component \( A \tox{\br{\id_A,\const_x}} A \times X \to Y \) is Borel, this follows immediately from the universal property of the product \( Y^X = \prod_{x \in X} Y \) and the fact that this product is preserved by \( U \colon \SBor \to \Set \).
\end{proof}

\section{The Borel Grothendieck topos}\label{sec:borel-topos}
We consider two pretopologies $J_{\SBor}^\inc \subset J_{\SBor}$ on $\SBor$.
The topology \( J_{\SBor} \), and the associated topos \( \Sh(J_{\SBor}) \), will be our main object of interest, but the topology \( J_{\SBor}^\inc \) is somewhat simpler, and many of our results holds in the topos \( \Sh(J_{\SBor}^\inc) \) as well, so we include it in our discussion.

Given a family \( \bA \in \Fam(X) \) of Borel maps, we declare that \( \bA \in J_{\SBor}(X) \) if \( \bA \) is countable, and the images of the $a \in \bA$ cover $X$.
The coverings in $J_{\SBor}^\inc$ are defined in the same way, except that the maps $a \in \bA$ in a covering \( \bA \) are also required to be injective.
In this case,
by the Lusin-Suslin theorem \cite[Theorem~15.1]{kechris-classical},
the image of each $a \in \bA$ is Borel; i.e., in $J_{\SBor}^\inc$, the coverings are countable coverings by Borel subsets.
One verifies easily that $J_{\SBor}^\inc$ and $J_{\SBor}$ are indeed pretopologies.
We also note that if, in the pretopology $J_{\SBor}^\inc$, we take only coverings by \emph{disjoint} Borel subsets---i.e., partitions---the resulting topology is the same.



\begin{defn}
  We define the topoi $\cE \subset \cE^{\inc} \subset \wh{\SBor}$ to be the categories $\Sh(\SBor,J_\SBor)$ and $\Sh(\SBor,J_\SBor^\inc)$ of sheaves on the presite $(\SBor,J_\SBor)$ and $(\SBor,J_\SBor^\inc)$, respectively.
\end{defn}

\begin{rmk}
  If we instead define the topology on \( \SBor \) so that the covering families are \emph{all} (not necessarily countable) jointly surjective families, the resulting topos is just equivalent to \( \Set \), see Proposition~\ref{propn:double-neg-subtopos}.
\end{rmk}

\begin{rmk} \label{rmk:coherent-topology}
  Each of the topologies \( J_{\SBor} \) and \( J_{\SBor}^\inc \) is an instance of a known general construction.
  The category \( \SBor \) is \emph{countably extensive}, and is in fact even the \emph{universal} Boolean countable complete countably extensive category (see \cite{chen-universal-characterization-of-standard-borel}).
  On any countably extensive category, there is a pretopology (the ``countable extensive pretopology'') in which the covering families are the countable coproduct diagrams \cite[Appendix~A]{roberts-internal-cats-ana-loca}.
  In the case of \( \SBor \), this is just \( J_{\SBor}^\inc \).

  Next, a family \( \bA \in \Fam(X) \) of morphisms in a category \( \cC \) is called \emph{extremal epimorphic} if any monomorphism \( P \tto X \) through which each \( (U_a \tox{a} X) \in \bA \) factors is an isomorphism.
  If the countable extremal epimorphic families in \( \cC \) are closed under pullbacks, they form a topology, called the \emph{countably coherent pretopology}, and on \( \SBor \), this is the topology \( J_{\SBor} \)---and in fact, in this case also agrees with the \( \aleph_1 \)-precanonical topology \cite[Chapter~9]{makkai-reyes-focl}, i.e., the countable families of the canonical topology (see \S\ref{subsec:can-top}).

  Usually, one considers the coherent topology only in the case when the underlying category is itself a \emph{countably coherent category}, which \( \SBor \) is not, since it lacks images.
  However, in \S\ref{subsec:analytic} below, we show that \( \cE \) \emph{is} the category of sheaves for the countably coherent topology on a countably coherent category, namely the category of \emph{analytic} spaces.
\end{rmk}

\begin{propn}
  The topology $J_\SBor$, and hence also $J_{\SBor}^\inc$ is subcanonical, i.e., $\yon \colon \SBor \to \wh \SBor$ factors through $\cE \hto \wh\SBor$.
\end{propn}

\begin{proof}
  It suffices to show that given maps $f \colon X \to Y$ and $g \colon W\to X$ in $\SBor$ with $g$ surjective, if $f \circ g \colon W\to Y$ is Borel, then so is $f$.
  This is the same as saying that for any subset $U \subset X$, if $g\I(U)$ is Borel, then so is $U$.
  But in this case, $U$ and $X \setminus U$ are both analytic,
  so this follows from Suslin's theorem
  \cite[Theorem~14.11]{kechris-classical}.
\end{proof}

Let \( \cD \) denote either of \( \cE \) or \( \cE^\inc \); most of what follows will apply in either case, and we will indicate when this is not so.
We have an ``underlying set'' functor \( \Gamma =  \cD(\tm,-) \), which, as with any topos, has a left adjoint \( \Delta \colon \Set \to \cE \) (the ``discrete space'' functor) taking \( X \in \Set \) to the sheafification of the constant presheaf \( \und X \) on \( X \).
It follows from \cite[Example~C3.6.3~(d)]{johnstone-elephant-vols} (and the fact that the terminal object in \( \SBor \) has no non-trivial covers) that \( \Gamma \) also has a right adjoint \( \nabla \colon \Set \to \cE \) (the ``indiscrete space'' functor), taking \( X \in \Set \) to the sheaf whose value at \( A \in \SBor \) is the set of all functions \( A \to X \), and that both \( \Delta \) and \( \nabla \) are fully faithful.
For \( X \) countable, the discrete space \( \Delta X \) agrees with the representable functor on \( X \):
\comm{I thought for a moment that also \( \yon \cong \nabla \circ U \) on countable spaces, but this is wrong: for any natural isomorphism \( \alpha \colon \yon \to \nabla \circ U \) (defined on the countable spaces): (i) \( \alpha_\tm \) is an iso (since its domain and codomain are terminal), (ii) whence it follows using naturality of \( \alpha \) along each \( \tm \to Y \) that \( (\alpha_Y)_\tm \colon \SBor(\tm, Y) \to Y^\tm \) is the identity for all (countable) \( Y \), (iii) whence it follows using naturality of \( \alpha_Y \) along each \( \tm \to X \) that \( (\alpha_Y)_X \colon \SBor(X,Y) \to Y^X \) is the inclusion map for all \( X \in \SBor \); (iv) and this inclusion is not an isomorphism when \( Y \) has at least two points and \( X \) is uncountable.}
\begin{propn}
  Given any countable set \( X \), the sheafification of the constant presheaf \( \und{X} \in \wh \SBor \) is the representable sheaf \( \wh X \in \cD \), with the universal map \( \alpha \colon \und{X} \to \wh X \) defined by \( \alpha_A(x) = \const_x \in \SBor(A,X) \) for \( A \in \SBor \).
\end{propn}

\begin{proof}
  Given \( F \in \cD \) and a morphism \( f \colon \und{X} \to F \) in \( \wh{\SBor} \), we need to show there is a unique \( \bar f \colon \wh X \to F \) with \( \bar f_A \circ \alpha_A = f_A \) for all \( A \in \SBor \).
  Given \( u \in \wh X(A) = \SBor(A,X) \), we have a covering \( \set{U_x \mid x \in X} \) of \( A \) in \( \SBor \) where \( U_x = u\I(x) \), and \( \rstr{u}{U_x} = \alpha_{U_x}(x) \in \wh X(U_X) \).
  If we are to have \( \bar f_{U_x} \circ \alpha_{U_x} = f_{U_x} \), since \( F \) is a sheaf, we must take \( \bar f_A(u) \in F(A) \) to be the unique element \( v \) with \( \rstr{v}{U_x} = f_{U_x}(x) \) for each \( x \in X \).
  This proves uniqueness of \( \bar f \), and existence follows from verifying that \( \bar f \) thus defined is indeed a morphism of sheaves.
\end{proof}
\noindent
For \( X \) countable, we thus have that \( \wh X \cong \coprod_{x \in X} \tm \).
From \cite[p.~270]{mac-lane-moerdijk}, we conclude:
\begin{cor}\label{cor:nat-is-nat}
  \( \wh{\N} \in \cD \) is a natural numbers object.
\end{cor}

Apart from the countable spaces, there is up to isomorphism a single object in \( \SBor \), namely \( \two^\N \), to which we turn our attention next.
\begin{propn}\label{propn:exps-are-exps}
  Given \( X,Y \in \SBor \) countable, the morphism \( \wh{Y^X} \times \wh X \toix{\br{\wh \pi_0,\wh \pi_1}\I} \wh{Y^X \times X} \tox{\wh \ev} \wh Y \) exhibits \( \wh{Y^X} \) as an exponential object of \( Y \) by \( X \) in \( \wh{\SBor} \), hence (by \cite[Proposition III.6.1]{mac-lane-moerdijk}) in \( \cD \).
\end{propn}
The fact used here that \( \br{\wh \pi_0,\wh \pi_1} \colon \wh{Y^X \times X} \to \wh{Y^X} \times \wh X \) is an isomorphism follows from the Yoneda embedding \( \yon \colon \SBor \to \wh{\SBor} \) being limit-preserving.
Proposition~\ref{propn:exps-are-exps} follows from Proposition~\ref{propn:sbor-exponentials} and the following general fact:
\comm{This is a very basic and general fact, but I couldn't find a reference for it.}
\begin{lem}
  For any small category \( \cC \), the Yoneda embedding \( \yon \colon \cC \to \wh \cC \) preserves exponentials.
\end{lem}
This means that for any exponential object \( \ev \colon Y^X \times X \to Y \) in \( \cC \), \( \wh \ev \colon \wh{Y^X \times X} \to \wh Y \) is an exponential object in \( \wh \cC \).
\begin{proof}
  Fix an exponential object \( \ev \colon Y^X \times X \to Y \) in \( \cC \).
  We have a commutative diagram
  \[
    \begin{tikzcd}
      \cC(A, Y^X) \ar[r, "\times X"] \ar[d, "\yon"] & \cC(A \times X, Y^X \times X) \ar[r, "\ev \circ"] \ar[d, "\yon"] & \cC(A \times X, Y) \ar[d, "\yon"] \\
      \wh\cC(\wh A, \wh{Y^X}) \ar[r, "\wh{\times X}"] & \wh\cC(\wh{A \times X}, \wh{Y^X \times X}) \ar[r, "\wh{\ev} \circ"] & \wh \cC(\wh{A \times X}, \wh Y).
    \end{tikzcd}
  \]
  Since the composite \( \cC(A, Y^X) \to \cC(X \times A, Y) \) of the top row is an isomorphism by assumption, and the vertical morphisms are isomorphisms, it follows that the composite \( \wh\cC(\wh A, \wh{Y^X}) \to \cC(\wh{X \times A}, \wh Y) \) of the bottom row is an isomorphism.
  This shows that \( \wh \ev \colon \wh{Y^X \times X} \to \wh Y \) has the desired universal property with respect to \emph{representable} objects.

  By \loccit, we know there exists an exponential object \( \wh Y^{\wh X} \) of \( \wh Y \) by \( \wh X \) in \( \wh \cC \).
  What we just showed implies that that the morphism \( f \colon \wh{Y^X} \to \wh Y^{\wh X} \) induced by \( \wh \ev \colon \wh{Y^X \times X} \to \wh Y \) induces an isomorphism \( f \circ \colon \wh \cC(\wh A, \wh{Y^X}) \to \wh \cC(\wh A, \wh Y^{\wh X}) \) for all \( A \in \cC \)---i.e., by Yoneda's lemma, an isomorphism \( f_A \colon \wh{Y^X}(A) \to \wh Y^{\wh X}(A) \) for all \( A \in \cC \)---and is thus an isomorphism of presheaves.
\end{proof}

\subsection{Subsheaves in \alttext{\( \cE \)}{E} and \alttext{\( \cE^{\inc} \)}{Eⁱⁿᶜ}} \label{subsec:subsheaves-in-e}
We continue to write \( \cD \) for either of \( \cE \) or \( \cE^{\inc} \).

Given \( X \in \SBor \), note that any subset \( A \subset X \) determines a subsheaf \( \wh A \) consisting of those \( f \colon Y \to X \) in \( \SBor \) with image in \( A \); if \( A \) is Borel, this agrees with the representable functor \( \wh A \).
Conversely, any subsheaf \( R \subset X \) yields a subset \( \norm{R} \), the \defword{support} of \( R \) consisting of those points \( x \in X \) for which the inclusion \( \set{x} \hto X \) in \( R \set{x} \).
It is immediate that \( R \subset \wh{\ \norm{R}\ } \) for \( R \subset \wh X \), and that \( A = \norm{\wh A} \) for \( A \subset X \).

\begin{propn}\label{propn:complement-in-rep}
  The complement \( \neg R \) of a subsheaf \( R \subset \wh X \) of a representable is given by \( \wh{\norm{R}^{\cm}} \), where \( \norm{R}^\cm = X \setminus \norm{R} \).
\end{propn}

\begin{proof}
  The claim is that \( \norm{R}^{\cm} \) has trivial intersection with \( R \), and that it is the maximal subsheaf with this property.
  Both of these facts are obvious, where we note that the bottom element \( \bot \in \Sub(\wh X) \) is given by \( \bot \emptyset = \Hom(\emptyset,X) \) (a singleton set) and \( \bot Y = \emptyset \subset \Hom(Y,X) \) for \( Y \) non-empty.
\end{proof}

\begin{cor}\label{cor:double-negations}
  For any subsheaf \( R \subset \wh X \), we have \( \neg \neg R = \wh{\ \norm{R}\ } \).
\end{cor}

We call a subsheaf \( R \subset F \) \defword{decidable} if \( R \vee \neg R = F \).

\begin{propn}\label{propn:decidable-sub-is-borel}
  \( R \subset \wh X \) is decidable iff \( R = \wh{A} \) for some Borel subset \( A \subset X \) (in which case \( A = \norm{R} \)).
\end{propn}

\begin{proof}
  The disjunction \( R \vee \neg R \) is the least subsheaf of \( \wh X \) containing both \( R \) and \( \neg R \).

  If \( R = \wh{A} \), then \( \neg R = \wh{A^\cm} \) by Proposition~\ref{propn:complement-in-rep}.
  Then, if \( S \subset \wh X \) is a subsheaf containing \( R \) and \( \neg R \), for any \( f \colon Y \to X \), we have a Borel cover \( Y = f\I(A) \cup f\I(A^\cm) \), and the restriction of \( f \) to each of these subsets is in \( R\pbig{f\I(A)} \) and \( (\neg R)\pbig{f\I(A^\cm)} \), respectively.
  Hence they are in \( S\pbig{f\I(A)} \) and \( S\pbig{f\I(A^\cm)} \), respectively, and thus \( f \in S(Y) \) since \( S \) is a sub-sheaf, and we conclude that \( S = \wh X \), as desired.

  Conversely, suppose \( R \vee \neg R = \wh X \).
  Then in particular, \( R = \neg \neg R \), hence \( R = \wh{A} \) where \( A = \norm{R} \).
  It remains to see that \( A \) is Borel.

  The assumption \( R \vee \neg R = \wh X \) gives in particular \( \id_X \in (R \vee \neg R) (X) \).
  It follows from the Kripke-Joyal semantics (\S\ref{subsec:kj-semantics}) that there is a cover \( a_j \colon U_j \to X \) for \( j = 0,1 \) with \( a_0 = a_0 \circ \id_X \in R(U_0) \) and \( a_1 = a_1 \circ \id _X \in (\neg R)(U_1) \).
  This means that the images of \( a_0 \) and \( a_1 \) are contained in \( A \) and \( A^\cm \) respectively, and hence are equal to \( A \) and \( A^\cm \), since these images cover \( X \).
  In the case \( \cD = \cE^{\inc} \), we are done, since the covers in \( J^\inc \) only involve Borel inclusions.
  And in the case \( \cD = \cE \), we have that \( A \) and \( A^\cm \) are both analytic and hence Borel by Suslin's theorem \cite[Theorem~14.11]{kechris-classical}.
\end{proof}

The above relationship between subsets of \( X \) and subsheaves of \( \wh X \) can be refined using \( \sigma \)-ideals.
In \( \cE \), this will allow us to classify all subsheaves of \( \wh X \).

Given \( X \in \SBor \), by a \defword{\( \sigma \)-ideal of analytic sets on \( X \)}, we mean a non-empty downward-closed set of analytic subsets of \( X \) which is closed under countable unions.
Let us write \( \AI(X) \) for the set of \( \sigma \)-ideals of analytic subsets \( X \).
We then have factorizations
\begin{equation}\label{eq:subsheaf-factorization}
  \begin{tikzcd}
    \pow(X) \ar[r, "I_{-}", shift left] &
    \AI(X) \ar[r, shift left, "\wh{(-)}"] \ar[l, shift left, "\abs{-}"] &
    \Sub(\wh X) \ar[l, "\abs{-}", shift left]
  \end{tikzcd}
  \qquad\text{of the maps}\qquad
  \begin{tikzcd}
    \pow(X) \ar[r, "\wh{(-)}", shift left] &
    \Sub(\wh X) \ar[l, "\norm{-}", shift left]
  \end{tikzcd}
\end{equation}
such that \( \abs{I_U} = U \) for \( U \in \pow(X) \) and \( \abs{\wh I} = I \) for \( I \in \AI(X) \), and in particular the rightward maps are injective and the leftward maps are surjective.

The maps displayed above are defined as follows:
\begin{itemize}
\item Given \( U \in \pow(X) \), \( I_U \in \AI(X) \) consists of those analytic subsets \( V \) with \( V \subset U \); this is clearly a \( \sigma \)-ideal.
\item For \( R \in \Sub(\wh X) \), \( \abs{R} \) consists of those \( V \subset X \) which are the image of some \( f \colon Y \to X \) with \( f \in R (Y) \); to see that \( \abs{R} \) is downward closed, given \( f \colon Y \to X \) in \( R(Y) \), any analytic subset \( U \subset f(Y) \) is of the form \( g(Z) \) for some \( g \colon Z \to X \), and is hence also the image of the map from the fibre product \( Y \times_X Z \to X \), which factors through \( f \), and is thus in the sieve \( R \); to see that \( \abs{R} \) is closed under countable unions, given a countable family of Borel maps \( Y_\alpha \to X \), each in \( R(Y_\alpha) \), the union of their images is the image of the Borel map \( \coprod_{\alpha} Y_\alpha \to X \), which is in \( R(\coprod_{\alpha} Y_\alpha) \) by the subsheaf condition.
\item For \( I \in \AI(X) \), \( \abs{I} \) is the set of points in \( I \), and
\item \( \wh I \tto \wh X \) is the following subsheaf: \( \wh I(Y) \) for \( Y \in \SBor \) consists of those \( f \colon Y \to X \) with \( f(Y) \in I \); this is a sieve by downward-closedness of \( I \), and is a subsheaf by the closedness of \( I \) under countable unions.
\end{itemize}

That these do indeed provide a factorization as in \eqref{eq:subsheaf-factorization} is easily seen.
It is also clear that \( \abs{I_U} = U \) for \( U \in \pow(X) \), as claimed.
That \( \abs{\wh{I}} = I \) for \( U \in\AI(X) \) is due to the fact that, by definition, every analytic set in \( X \) is the image of some map \( Y \to X \) in \( \SBor \).

\begin{propn}\label{propn:subsheaves-in-e}
  For \( X \in \SBor \), the above map \( \wh{(-)} \colon \AI(X) \to \Sub_{\cE}(\wh X) \) is surjective (hence bijective by what was already said), i.e., every subsheaf of \( \wh X \) is of the form \( \wh I \) for a unique \( \sigma \)-ideal \( I \in \AI(X) \).
\end{propn}

\begin{proof}
  Given any subsheaf \( R \subset \wh X \), we need to show that \( R = \wh{\abs{R}} \).
  The inclusion \( R \subset \wh{\abs{R}} \) is immediate.
  Hence given any \( f \colon Y \to X \) with \( f(Y) \in \abs{R} \), we need to show that \( f \in R(Y) \).
  By assumption, there is some \( g \colon Z \to X \) in \( R \) with \( f(Y) = g(Z) \).
  It follows that the map from the fibre product \( Y \times_X Z \to X \) is in \( R \).
  Moreover, it factors as \( Y \times_X Z \tox{\pi_Y} Y \tox{f} X \), where \( \pi_Y \) is surjective and thus a cover in \( J_{\SBor} \).
  Since \( R \) is a subsheaf, we conclude that \( f \in R(Y) \), as desired.
\end{proof}

\begin{rmk}\label{rmk:subsheaf-counterexamples}\enumbelow
  \begin{enumerate}
  \item For \( X \) countable, the map \( I_- \colon \pow(X) \to \AI(X) \) is also surjective, as is easily seen.
    However, if \( X \) is uncountable, then it is not, and thus the subsheaves of \( \wh X \) in \( \cE \) are not exhausted by those of the form \( \wh{I_A} \) for \( A \subset X \).
    An example is the \( \sigma \)-ideal consisting of the countable subsets of \( X \).
  \item Proposition~\ref{propn:subsheaves-in-e} fails in \( \cE^\inc \).
    An example is obtained by taking a surjective map \( p \colon X' \to X \) in \( \SBor \) which does not admit a section (see \cite[Exercise~18.17]{kechris-classical}),
    and taking \( R \subset \wh X \) to be the subsheaf consisting of all maps \( f \colon Y \to X \) which factor through \( p \).
    Then \( R \) is certainly a sieve, and it is a subsheaf, since given a map \( f \colon Y \to X \), if \( f \) admits a lift to \( X' \) after restricting to each set in a Borel partition of \( Y \), then these can be assembled to give a lift \( Y \to X' \).
    However, \( R \ne \wh{\abs{R}} \), since if these were equal, then \( \id_X \colon X \to X \) would lift to \( X' \), which by assumption it does not.
  \end{enumerate}
\end{rmk}

\subsection{Analytic spaces} \label{subsec:analytic}
We now give an alternate description of \( \cE \) as sheaves on the category of \emph{analytic spaces}, as promised in Remark~\ref{rmk:coherent-topology}.
\begin{defn}
  An \defword{analytic space} is a measurable space given by a subset \( A \) of a complete Borel space \( X \) with the induced \( \sigma \)-algebra (consisting of sets \( A \cap B \subset A \) with \( B \subset X \) Borel).
  We denote by \( \SAn \) the full subcategory of \( \Meas \) spanned by the analytic spaces; thus, \( \SBor \subset \SAn \subset \Meas \).

  Given \( A \subset X \) with \( A \in \SAn \) an analytic subset of \( X \in \SBor \), we call a subset \( V \subset A \) \defword{analytic} if it is analytic in \( X \), or equivalently, if it is the image of some \( B \to A \) in \( \SAn \).
\end{defn}
\noindent
The morphisms in \( \SAn \) can be described in several equivalent ways:
\begin{propn} \label{propn:san-mor-chars}
  Let \( X,Y \in \SBor \) and let \( i_A \colon A \hto X \) and \( i_B \colon B \hto Y \) be analytic subsets.
  For a map \( f \colon A \to B \), the following are equivalent:
  \begin{enumerate}[(i)]
  \item \( f \) is a morphism in \( \SAn \) (i.e., the preimage under \( f \) of any Borel subset of \( B \) is Borel)
  \item The preimage under \( f \) of any analytic set is analytic.
  \item \label{item:ana-as-subsheaf} For any \( u \colon Z \to X \) in \( \SBor \) factoring through \( A \), the map \( i_B \circ f \circ u \colon Z \to Y \) is Borel.
  \item For \emph{some} \( u \colon Z \to X \) in \( \SBor \) with image \( A \), the map \( i_B \circ f \circ u \colon Z \to Y \) is Borel.
  \end{enumerate}
\end{propn}
\noindent

\begin{proof}
This can be proven, e.g., in the order
(iii) \( \To \) (iv) \( \To \) (ii) \( \To \) (i) \( \To \) (iii).
The only non-trivial step is
The implication (ii) \( \To \) (i), which follows from Lusin's separation theorem \cite[Theorem~14.7]{kechris-classical}.
\end{proof}

The characterization \ref{item:ana-as-subsheaf} shows that the morphisms \( A \to B \) are in bijection with the morphisms \( \wh{I_A} \to \wh{I_B} \) between the subsheaves \( \wh{I_A} \subset \wh X \) and \( \wh{I_B} \subset \wh Y \) of \S\ref{subsec:subsheaves-in-e}.
It follows that \( \SAn \) can be alternatively characterized as the closure under images of the image of the embedding \( \yon \colon \SBor \to \cE \):
\begin{propn}
  Define \( F \colon \SAn \to \cE \) by \( F(A) = \wh{I_A} \) for \( i_A \colon A \hto X \in \SBor \) analytic, and for \( f \colon A \to B \), define \( F(f) \) by mapping \( u \in \wh{I_A}(Z) \subset \SBor(Z, X) \) to \( i_B \circ f \circ u \).
  Then \( F \) is a fully faithful functor extending \( \yon \colon \SBor \to \cE \) whose essential image consists of the images of morphisms between representables in \( \cE \).
\end{propn}
\begin{proof}
  The functoriality of \( F \) is obvious, as is the fact that it extends \( \yon \).
  The only thing we need to check is that if \( A \subset X \) is the image of \( h \colon Y \to X \), then \( \wh{I_A} \subset \wh X \) is the image of \( \hat h \colon \wh X \to \wh Y \) in \( \cE \), i.e., that \( F(h) \colon \wh Y \to \wh{I_A} \) is an epi and \( F(i_A) \colon \wh{I_A} \to \wh X \) is a mono.
  The latter is clear, since this is just the inclusion \( \wh{I_A} \hto \wh X \).

  For the former, given \( U \in \SBor \) and \( k \colon \wh U \to \wh{I_A} \), we need to find an epi \( a \colon U' \too U \) and \( \hat \ell \colon \wh{U'} \to \wh Y \) making the following square commute.
  \[
    \begin{tikzcd}
      \wh{U'} \ar[r, "\hat \ell"] \ar[d, "\hat a"] & \wh Y \ar[d, "F(h)"] \\
      \wh U \ar[r, "k"] & \wh{I_A}
    \end{tikzcd}
  \]
  We may take \( U' \) to be the fibre product \( U' = U \times_X Y \tox{a = \pi_U} U \) and \( \ell = \pi_Y \colon U' \to Y \).
\end{proof}
\begin{rmk}
  The subcategory of \( \cE \) described above is (up to equivalence) indeed the \emph{closure under images} of \( \yon(\SBor) \)---i.e., the least subcatgeory of \( \cE \) closed under images containing \( \yon(\SBor) \).
  That is to say, to close under images, it is enough to take images of all morphisms \emph{once}: given images \( \im(f) \) and \( \im(g) \) of \( f \colon P \to Q \) and \( g \colon R \to S \), the image of any morphism \( \im(f) \to \im(g) \) is also the image of a morphism \( P \to S \).
\end{rmk}

We next want to show that \( \SAn \) has countable limits and well-behaved images and countable coproducts, and (hence) also well-behaved countable joins of subobjects.
Here, that images and joins are well-behaved just means they are stable under pullbacks, and for coproducts, it means that they are stable under pullbacks and are ``disjoint'' (see, e.g., \cite{chen-universal-characterization-of-standard-borel}).
There are a variety of names for various combinations of these conditions: a category \( \cC \) is
\begin{itemize}
\item \emph{regular} if it is has finite limits and stable images
\item \emph{coherent} (or \emph{logical} in \cite{makkai-reyes-focl}) if it is regular and has stable finite joins
\item \emph{countably coherent} (or \emph{\( \omega_1 \)-logical} in \cite{makkai-reyes-focl} and \emph{\( \sigma \)-coherent} in \cite{johnstone-elephant-vols}) if it is regular and has stable countable joins
\item (countably) \emph{extensive} if has stable disjoint finite (or countable) coproducts.
\end{itemize}
Sometimes (e.g., in \cite{roberts-internal-cats-ana-loca}) ``lextensive'' is also used for extensive and finitely complete, and in \cite{johnstone-elephant-vols}, ``(\( \sigma \)-)positive (\( \sigma \)-)coherent'' is used for regular and (\( \sigma \)-)extensive (which implies (\( \sigma \)-)coherent).

We also need one more concept: we call a category \defword{regularly countably complete} if it is regular and countably complete, and any countable product of regular epis is again a regular epi (and hence the image of a product of morphisms is the product of the images).

\begin{propn} \label{propn:san-is-good}
  The category \( \SAn \) is regularly countably complete and countably extensive (and hence countably coherent), and the underlying set functor \( \SAn \to \Set \) is regularly countably complete and countably extensive (in the sense described below); moreover the epis in \( \SAn \) (and in \( \SBor) \) are precisely the surjections, and are all regular.
\end{propn}

\begin{proof}
  That \( \SAn \) is regular follows from its being the closure under images in \( \cE \) of the finite-limit category \( \SBor \).
  Next, that \( \SAn \) has stable countable coproducts similarly follows since \( \SBor \) has countable coproducts, and both monos and epis, and hence images, are closed under arbitrary coproducts in any Grothendieck topos; moreover, the coproducts are disjoint since they are in \( \cE \).
  To see that \( \SAn \) has countable products, given analytic subsets \( A_i \subset X_i \in \SBor \), we have the analytic subset \( \prod_i A_i \subset \prod_i X_i \) (the image of \( \prod_i f_i \), where \( A_i \) is the image of \( f_i \colon Y_i \to X_i \) with \( Y_i \in \SBor \)), and it follows immediately from Proposition~\ref{propn:san-mor-chars} that this has the universal property of the product.
  That regular epis are closed under countable products will follow immediately once we show that these are precisely the surjections and that \( U \colon \SAn \to \Set \) preserves products.

  Now, that \( U \) preserves limits is immediate from the fact that it is isomorphic to the hom functor \( \SAn(\tm, -) \) from the terminal object, and that it preserves countable coproducts follows from this and the facts that the coproducts are disjoint and that \( \tm \) has no non-trivial coproduct decompositions.
  That \( U \) preserves images will follow once we show that the regular epis are the surjections.

  That every surjection is an epi is obvious; conversely, if \( f \colon X \to Y \) is not surjective, take \( y \in Y - \im(f) \) and \( g_0,g_1 \colon Y \to \two \) with \( g_0(Y - \set{y}) = g_1(Y - \set{y}) = 0 \) and \( g_i(y) = i \); then \( g_0 f = g_1 f \) but \( g_0 \ne g_1 \).
  It remains to see that every epi is regular, i.e., that if \( f \colon X \to Y \) is epi, then \( X \times_Y X \toto X \to Y \) is a coequalizer.
  In \( \SBor \), this is an instance of the subcanonicity of the topology \( J_{\SBor} \), which we verified in \S\ref{sec:borel-topos}, and the verification in \( \SAn \) is similar, using Proposition~\ref{propn:san-mor-chars}.
\end{proof}

Next, we show that \( \SAn \) is the \emph{completion} of \( \SBor \) in a suitable sense, and in fact in a number of different senses.
It is the \emph{regular} completion, and also the \emph{regularly countably complete, countably extensive} completion.
By the first of these, we mean that \( \SAn \) is regular, that the inclusion \( i \colon \SBor \hto \SAn \) is regular (i.e., preserves finite limits and whatever images exist in \( \SAn \)), and that for any other regular category \( \cC \), precomposition with \( i \) induces an equivalence
\begin{equation} \label{eq:reg-comparison}
  \Fun^\reg(\SAn, \cC) \to \Fun^\reg(\SBor, \cC)
\end{equation}
between categories of regular functors.
By the second, we mean the same thing, but replacing ``regular'' with ``regular, countably complete, countably extensive''.
Here, by a ``countably complete'' resp.\ ``(countably) extensive functor'' functor, we just mean one preserving countable limits resp.\ (countable) coproducts.
Similarly, a ``(countably) coherent'' functor is a regular functor preserving (countable) joins of subobjects.

And in fact, we will see that \( \SAn \) is also the completion of \( \SBor \) with respect to any of the several intermediate notions between ``regular'' and ``regularly countably complete, countably extensive'', namely any combination of ``regularly countably complete'', ``coherent'', ``countably coherent'', ``regular and extensive'', and ``regular and countably extensive''.
We note the phrase ``completion of a category \( \cC \) under \( X \)'s'' often refers to \emph{freely} adjoining \( X \)'s to \( \cC \) (while possibly preserving some other structure in \( \cC \)), whereas here, we always mean adjoining \( X \)'s to \( \cC \) while \emph{preserving any \( X \)'s already present in \( \cC \)}.
Nonetheless, the \emph{construction} of the completions we give---i.e., taking a closure inside of a sheaf category---is a common one, and is used, e.g., for the exact completion of a regular category in \cite{lack-exact-completion}.

\begin{thm} \label{thm:reg-completion}
  \( \SAn \) is the regular completion of \( \SBor \), the regularly countably complete, countably extensive completion of \( \SBor \), and the completion of \( \SBor \) in all of the intermediate senses listed above.
\end{thm}
\comm{It almost seems that this should imply a universal property of \( \SAn \), using the universal property of \( \SBor \) from \cite{chen-universal-characterization-of-standard-borel}, but it doesn't seem to work, since the latter involves \emph{boolean} categories, and \( \SAn \) is not boolean.}
\noindent
See Appendix~\ref{sec:reg-completion-proof} for the proof.

\begin{rmk}
  Passing from \( \SBor \) to its coherent completion \( \SAn \) amounts to allowing the interpretation of existential quantifiers, alongside the conjunctions and disjunctions which were already available in \( \SBor \) (note that negations were also available in \( \SBor \), but this is no longer so in \( \SAn \)).
  Since the topos \( \cE \) is of course coherent, we obtain an induced functor \( \SAn \to \cE \), and as we observed above, this is fully faithful, enabling us to identify the coherent completion with a full subcategory of \( \cE \).
  Now in \( \cE \), we can also interpret the rest of first-order logic, i.e., universal quantifiers and implications.
  We can correspondingly consider the \emph{Heyting} completion \( \cC \) of \( \SBor \), which will come with a functor \( \cC \to \cE \).
  One can again ask if this is fully faithful, so that the Heyting completion can be recovered as the closure of \( \SBor \) in \( \cE \) under the Heyting operations; we do not know whether this is true.
\end{rmk}

We end this section by observing, as promised in Remark~\ref{rmk:coherent-topology}, that \( \cE \) is the equivalent to the topos of sheaves for the countably coherent topology on the countably coherent category \( \SAn \); this is just the pretopology \( J_{\SAn} \) in which the covers are, again, the countably jointly surjective families (and, as with \( J_{\SBor} \), it is also the \( \omega_1 \)-precanonical topology).
More precisely, since the restriction of \( J_{\SAn} \) to \( \SBor \) is \( J_{\SBor} \), every sheaf on \( (\SAn, J_{\SAn}) \) restricts to one on \( (\SBor, J_{\SBor}) \) giving a restriction functor \( \Sh(\SAn, J_{\SAn}) \to \cE \).

\begin{propn}
  The restriction functor \( \Sh(\SAn, J_{\SAn}) \to \cE \) is an equivalence.
\end{propn}

\begin{proof}
  This follows immediately from the ``comparison lemma'' for dense sub-sites \cite[Theorem~III.4.1]{sga-4-1}.
\end{proof}

\subsection{Other topos-theoretic properties} \label{subsec:more-topos-props}
We gather some other facts about the topoi \( \cE \) and \( \cE^\inc \).
Again let \( \cD \) denote either one.

\begin{propn}\label{propn:two-valued}
  \( \cD \) is two-valued, i.e., the terminal object has precisely two subobjects.
\end{propn}

\begin{proof}
  The terminal object in \( \cD \) is represented by the one-point space \( \tm \) in \( \SBor \), and in fact there are only two sieves \( R \subset \wh{\tm} \) (both of which are subsheaves).
\end{proof}

\begin{propn}
  \( \cD \) is not boolean (i.e., not every subobject lattice in \( \cD \) is a Boolean algebra).
\end{propn}

\begin{proof}
  This follows from Remark~\ref{rmk:subsheaf-counterexamples}, which gives examples of undecidable subobjects.
\end{proof}

\begin{cor}\enumbelow
  \begin{enumerate}[(i)]
  \item \( \cD \) does not satisfy the axiom of choice (i.e., not every epimorphism in \( \cD \) has a section).
  \item \( \cD \) is not well-pointed (i.e., there exist distinct parallel morphisms \( f,g \colon X \to Y \) such that \( f \circ a = g \circ a \) for every \( a \colon \tm \to X \)).
  \end{enumerate}
\end{cor}

\begin{proof}
  These follow from \cite[Theorem~5.23~and~Proposition~9.33]{johnstone-topos-theory}, respectively.
\end{proof}

\begin{rmk}
    Despite the above proposition, a certain limited form of choice does hold in \( \cD \); see \S\ref{subsubsec:countable-choice}.
\end{rmk}

\begin{propn}
  Markov's principle \cite[\S8]{hyland-effective-topos} holds in \( \cE \); that is:
  \[
    \cE \vDash \forall P \tolon \pow(\wh \N).\
    (\forall n \tolon \N.\ n \in P \vee \neg n \in P) \wedge
    \neg \neg \exists n \tolon \N.\ n \in P \to
    \exists n \tolon \N.\ n \in P.
  \]
\end{propn}
\noindent
Here, we recall that the power sheaf \( \pow(\wh \N) \) (which is a power object in the sense of \S\ref{subsec:elementary-topoi}) is the sheaf taking \( U \in \SBor \) to the set \( \Sub(\wh{U} \times \wh{\N}) \) of subsheaves of \( \wh{U} \times \wh{\N} \), and taking a morphism \( f \colon U \to V \) to the pre-image map \( (f \times \id)\I \colon \Sub(\wh V \times \wh \N) \to \Sub(\wh U \times \wh \N) \).
The subsheaf \( \unin \hto \wh \N \times \pow(\wh \N) \) takes \( U \in \SBor \) to the subset of \( \SBor(U,\N) \times \Sub(\wh U \times \wh \N) \) consisting of those pairs \( (f, R) \) such that \( \br{\id, \hat f} \colon \wh U \to \wh U \times \wh \N \) factors through \( R \hto \wh U \times \wh N \).
\begin{proof}
  We must show that for any \( U \in \SBor \), any \( P \colon U \to \pow(\N) \) satisfying
  \begin{equation} \label{eq:markov-decide}
    U \Vdash \forall n \tolon \N.\ n \in P \vee \neg n \in P
\end{equation}
  and
  \begin{equation} \label{eq:markov-doubleneg}
    U \Vdash \neg \neg \exists n \tolon \N.\ n \in P
  \end{equation}
  also satisfies
  \begin{equation} \label{eq:markov-concl}
    U \Vdash \exists n \tolon \N.\ n \in P.
  \end{equation}
  By Proposition~\ref{propn:decidable-sub-is-borel}, the condition \eqref{eq:markov-decide} is equivalent to the subsheaf of \( \wh U \times \wh \N \) classified by \( U \) being of the form \( \wh R \) for some Borel subset \( R \subset U \times \N \).
  Since every non-empty \( U \in \SBor \) admits a morphism from \( \tm \), and every \( U \to \tm \) with \( U \) non-empty is a cover, it follows that condition \eqref{eq:markov-doubleneg} is equivalent to there being, for each \( u \colon \tm \to U \), some \( n \colon \tm \to \N \) such that \( \br{u,n} \colon \tm \to U \times \N \) lies in \( R \)---in other words, to the composite \( R \htox{i_R} U \times \N \tox{\pi_U} \) being surjective.
  Finally, the condition \eqref{eq:markov-concl} we wish to prove says that there exists some cover of \( U \)---i.e., some Borel surjection \( p \colon U' \to U \)---such that, letting \( R' = (p \times \id_\N)\I(R) \subset U' \times \N \), the composite \( R' \htox{i_{R'}} U' \times \N \tox{\pi_{U'}} U' \) has a section.
  But since we are assuming that \( \pi_U \circ i_R \) is surjective, we can simply take \( U' = R \) and \( p = \pi_U \circ i_R \).
  Then \( R' = \set{((u,m),n) \in R \times \N \mid (u,n) \in R} \), and thus \( \pi_{U'} \circ i_{R'} \colon R' \to R \) has a Borel section \( R \tox{\br{\id_R,\pi_\N \circ i_R}} R' \subset R \times \N \).
\end{proof}

\begin{rmk}
  We suspect that Markov's principle does \emph{not} hold in \( \cE^\inc \), but we have not verified this.
\end{rmk}

\subsection{The double negation sub-topos} \label{subsec:double-neg-subtopos}
Fix a Grothendieck topos \( \cD \).
Recall (\cite[\S V.1]{mac-lane-moerdijk}) that an object \( Y \in \cD \) is called a \emph{double-negation sheaf} if for every \( X \in \cD \) and every subsheaf \( R \subset X \) with \( \neg \neg R = X \), composition with the inclusion \( i_R \colon R \hto Y \) induces a bijection \( i_R^* \cD(X, Y) \toi \cD(R, Y) \).
The full subcategory \( \cD_{\neg \neg} \) consisting of double-negation sheaves is itself a Grothendieck topos, and the inclusion \( \cD_{\neg \neg} \hto \cD \) has a left adjoint \( \cD \to \cD_{\neg \neg} \) (``double-negation-sheafification'') preserving finite limits.
\begin{lem}
  An object \( Y \in \cD \) is a double-negation sheaf if and only if it satisfies the above condition \( \neg \neg R = X \To \cD(X, Y) \toi \cD(R, Y) \) for all subsheaves \( R \subset X \) of \emph{representable} sheaves \( X \in \cD \).
  \comm{This must be well-known, but I was unable to find a direct statement of it in \cite{mac-lane-moerdijk,johnstone-elephant-vols,johnstone-topos-theory}.}
  \comm{This lemma has nothing to do with double-negation sheaves; it's true for any Lawvere-Tierney topology.}
\end{lem}

\begin{proof}
  One direction is trivial.
  The other direction follows by observing that the class of objects \( X \) for which this condition holds is closed under colimits, and that every \( X \in \cD \) is a colimit of representables.
\end{proof}

We now specialize again to \( \cD \) being of \( \cE \) or \( \cE^\inc \).
On \( \SBor \), consider the pretopology \( J_{\SBor}^{\neg \neg} \) in which the covers are arbitrary (not necessarily countable) jointly surjective families.
Thus \( J_{\SBor}^\inc \subset J_{\SBor} \subset J_{\SBor}^{\neg \neg} \), and hence \( \Sh(\SBor, J_{\SBor}^{\neg \neg}) \subset \cE \subset \cE^{\inc} \).
\begin{propn}\label{propn:double-neg-subtopos}
  In either of \( \cD = \cE \) or \( \cD = \cE^\inc \), the double negation sub-topos \( \cE_{\neg \neg} \) is equal to \( \Sh(\SBor, J_{\SBor}^{\neg \neg}) \), and also to the essential image of the ``indiscrete space'' functor \( \nabla \colon \Set \to \SBor \), and is thus equivalent to \( \Set \).
\end{propn}

\begin{proof}
  We claim that the topology on \( \SBor \) generated by the pretopology \( J_{\SBor}^{\neg \neg} \) is also generated by the collection \( K \) of covering families \( \bA \in \Fam(X) \) consisting of the inclusions of all the points in \( X \) (i.e., each \( U_a \) for \( a \in \bA \) is a singleton and the images of the \( a \colon U_a \to X \) are all disjoint).
  Indeed, \( K \subset J_{\SBor}^{\neg \neg} \), and by axiom (T4) from \S\ref{subsection:topologies}, \( J_{\SBor}^{\neg \neg} \subset \br{K} \), and hence \( \br{K} = \br{J_{\SBor}^{\neg \neg}} \).

  Hence, since \( K \) satisfies axiom (T2), \( \Sh(\SBor, J_{\SBor}^{\neg \neg}) \) consists of those presheaves \( F \colon \SBor^\op \to \Set \) which are sheaves with respect to \( K \).
  But any such presheaf is clearly determined by \( F(\tm) \), since \( F(X) \cong \prod_{x \in X} F(\tm) \) by the sheaf condition, and is thus isomorphic to \( \nabla\pbig{F(\tm)} \)---and conversely, \( \nabla Z \) is clearly a \( K \)-sheaf for any \( Z \in \Set \).
  Thus, \( \Sh(\SBor, J_{\SBor}^{\neg \neg}) \) is the essential image of \( \nabla \).

  It remains to see that \( Y \in \cD \) is a double-negation sheaf iff it is a sheaf for \( J_{\SBor}^{\neg \neg} \).
  In general, any family of morphisms \( \bA \in \Fam(X) \) generates a sieve \( i_R \colon R \hto \wh X \) (consisting of those \( Y \to X \) factoring through some \( a \colon U_a \to X \) in \( \bA \)), and a presheaf \( F \in \wh{\SBor} \) is a sheaf with respect to \( \bA \) if and only if \( i_R^* \colon \wh{\SBor}(\wh X, F) \to \wh{\SBor}(R, F) \) is a bijection.
  By the universal property of sheafification, it follows that \( F \in \cD \) is a sheaf with respect to \( \bA \) iff \( i_{R'}^* \colon \wh{\SBor}(\wh X, F) \to \wh{\SBor}(R', F) \) is a bijection, where \( R' \in \cD \) is the sheafification of \( R \).

  By Corollary~\ref{cor:double-negations}, a subsheaf \( R \subset \wh X \) of a representable satisfies \( \neg \neg R = X \) if and only if each point \( \tm \to X \) is in \( R \), which is equivalent to being the sheafification of the sieve generated by some jointly surjective family. \comm{This is because sheafification increases sieves: \( R \subset R' \).}
  Thus, \( Y \in \cD \) is a double-negation sheaf iff it is a sheaf for \( J_{\SBor}^{\neg \neg} \), as required.
\end{proof}

\subsection{The canonical topology on \alttext{\( \SBor \)}{SBor}} \label{subsec:can-top}
Another natural topology to consider on \( \SBor \) is the \emph{canonical} topology, which is by definition the maximal subcanonical topology, i.e., the maximal topology for which every representable functor is a sheaf \cite[p.~15]{johnstone-topos-theory}.

In this section, we prove:
\begin{thm}\label{thm:canonical}
  Assuming the continuum hypothesis, the canonical topology is equal to \( J_{\SBor} \).
\end{thm}
We will show below that the assumption of the continuum hypothesis is necessary here, in the sense that it is consistent with \( \neg \mathrm{CH} \) that \( J_{\SBor} \) is not the canonical topology.

Recall that a family \( \bA \in \Fam(X) \) is \emph{effective epimorphic} if each representable \( \wh Y \) is a sheaf with respect to \( \bA \).
The canonical topology consists exactly of those families \( \bA \in \Fam(X) \) which are effective epimorphic and such that the pullback family \( f^* \bA \in \Fam(X') \) is effective epimorphic for every \( f \colon X' \to X \) in \( \SBor \).

Given \( X,Y \in \SBor \) and an analytic subset \( U \subset X \), recall from \S\ref{subsec:analytic} that a function \( f \colon U \to Y \) is \emph{Borel} if the preimage of each Borel set under \( f \) is analytic---equivalently, if \( f \circ g \colon Z \to Y \) is Borel for any \( Z \in \SBor \) and Borel \( g \colon Z \to X \) with image \( U \).
We have the following characterization of effective epimorphic families in \( \SBor \), which follows readily from the definitions.
\begin{propn}
  Fix a family \( \bA \in \Fam(X) \) in \( \SBor \), and consider the resulting family \( \cU = \set{a(U_a) \mid a \in \bA} \) of analytic subsets of \( X \).
  Then \( \bA \) is effective epimorphic if and only if \( \bigcup \cU = X \), and given any function \( f \colon X \to Y \) with \( Y \in \SBor \), if the restriction of \( f \) to each \( U \in \cU \) is Borel, then \( f \) is Borel.
\end{propn}

Theorem~\ref{thm:canonical} is thus reduced to:
\begin{propn}
  (Assuming \textup{CH}.)
  Given \( X \in \SBor \) and an uncountable cover \( \cU \) of \( X \) by analytic sets with no countable subcover, there is a Borel map \( g \colon X' \to X \) and a function \( f \colon X' \to Y \) with \( Y \in \SBor \) which is not Borel but whose restriction to \( g\I(U) \) is Borel for each \( U \in \cU \).

  In fact, we can take \( g = \id_X \) and $Y = \set{0, 1}$.
\end{propn}
\begin{proof}
  We will build a set $W \subseteq X$ and take $f$ to be the indicator function of $W$.
  We will ensure that $W$ is not Borel, but that for each $U \in \cU$, $W\cap U$ is countable and hence Borel.
  This is enough to ensure that $f$ is not Borel but that its restriction to each $U \in \cU$ is Borel.

  We will build $W$ by a fairly standard transfinite recursion of length $\omega_1$.
  Note that there are continuum many Borel subsets of $X$ as well as continuum many analytic subsets of $X$.
  Hence by CH we may fix length $\omega_1$ enumerations $\set{U_\alpha}_\alpha$ and $\set{V_\alpha}_\alpha$ of the elements of $\cU$ and of the Borel subsets of $X$, respectively (note that since $\cU$ is a cover of $X$ with no countable subcover, we must have $\# \cU = \omega_1$, though this does not actually matter).

  On each step of our transfinite recursion, we will put at most one point into $W$ and pick at most one point which we commit to never adding to $W$ on a future step.
  Aside from the points we have decided to commit to never adding to $W$, we will also require that for each $\alpha < \omega_1$, we do not add any point in $U_\alpha$ to $W$ on any step after $\alpha$.
  Since at most countably many points could have been added to $W$ by step $\alpha$ (at most one for each $\beta < \alpha$), this will ensure that $W \cap U_\alpha$ is countable.

  We will now describe what happens at step $\alpha$ of the transfinite recursion.
  Recall that we need $W$ to not be Borel.
  To ensure this, on step $\alpha$ we will take action to ensure $W \neq V_\alpha$.
  There are two cases to consider, depending on whether $V_\alpha$ is countable or not.

  First, suppose that $V_\alpha$ is uncountable.
  Note that we have only added countably many points to $W$ so far.
  Thus we can pick some point in $V_\alpha \setminus W$ and commit to never add it to $W$ in the future.
  This ensures that $W \neq V_\alpha$.

  Next, suppose that $V_\alpha$ is countable.
  Let $Z$ denote the countably many elements which we have already committed to not adding to $W$.
  We claim that there is some point in $X$ which is not in $Z\cup V_\alpha \cup \bigcup_{\beta < \alpha}U_\beta$---in other words, some point in $X \setminus V_\alpha$ which it is valid for us to add to $W$ on this step.
  If not, then $Z \cup V_\alpha \cup \bigcup_{\beta < \alpha}U_\beta$ covers $X$.
  Since $Z \cup V_\alpha$ is countable, we can pick countably many element of $\cU$ which cover it and which, together with $\set{U_\beta}_{\beta < \alpha}$, constitute a countable cover of $X$ by elements of $\cU$.
  But this contradicts our assumption about $\cU$.
  Thus we can pick some point in $X$ which is not in $V_\alpha$ and add it to $W$, which ensures that $W \neq V_\alpha$.

  This finishes our description of $W$.
  Our construction of $W$ ensures that it is not equal to any $V_\alpha$ and hence not Borel.
  And as we already explained, for each $U \in \cU$, $W\cap U$ is countable.
\end{proof}

We next address the necessity of CH above:
\begin{propn}[Hou]
\label{propn:not_canonical}
  The following is consistent with ZFC:
  there exists an uncountable cover \( \cU \) of \( \R \) by analytic sets with no countable subcover such that, given any Borel map \( g \colon X \to \R \) and any function \( f \colon X \to Y \) with \( X,Y \in \SBor \), if the restriction of \( f \) to \( g\I(U) \) for each \( U \in \cU \) is Borel, then so is \( f \).
\end{propn}
The proof of this proposition, which is due to Edward Hou, can be found in Appendix~\ref{sec:necessity-of-ch}; see \cite{hou-mo-forcing-argument} for a related argument.

\subsection{Van Osdol's topos}\label{subsec:van-osdol-topos}
In \cite{van-osdol-borel-topos}, a different ``Borel topos'', is considered.
It is the category \( \Sh(\cS, J_{\cS}) \) of sheaves on the non-full category \( \cS \subset \SBor \) consisting of Borel subsets of the unit interval \( [0,1] \), and inclusion functions between these, where \( J_{\cS} \) is the topology given by countable coverings.

The problem with this topos, from our point of view, is that---contrary to what is claimed in \opcit---the natural functor \( \Phi \colon \SBor \to \Sh(\cS, J_{\cS}) \) taking \( X \in \SBor \) to \( \SBor(-,X) \colon \cS^\op \to \Set \) is not fully faithful; specifically, it is not full.
To see why, let \( X, Y \in \SBor \), and fix a Borel function \( F \colon [0, 1] \times X \to Y \).
We can then define a natural transformation \( \alpha \colon \Phi X \to \Phi Y \) by defining \( \alpha_U \colon \SBor(U,X) \to \SBor(U,Y) \) for each Borel \( U \subseteq [0, 1] \) to be the map taking \( g \colon U \to X \) to the composite \( U \tox{\br{\id, g}} [0,1] \times X \tox{F} Y \).
However, this is only in the image of \( \Phi \) if \( F \) factors through the projection \( [0,1] \times X \to X \)---as can be seen, for instance, by taking \( U = [0,1] \) and \( g = \const_x \) where \( x \in X \) satisfies \( F(t,x) \ne F(t',x) \) for some \( t,t' \in [0,1] \).
\comm{It might be interesting to understand how \( \Sh(\cS, J_{\cS}) \) is related to our topos; in particular, the inclusion \( i \colon \cS \hto \SBor \) is a morphism of sites (with respect to \( J_\SBor^\inc \) and hence also \( J_\SBor \)), and hence induces a geometric morphism \( i_* \colon \cE^\inc \to \Sh(\cS, J_{\cS}) \) (which is not, however, fully faithful).}

\section{Borel mathematics and the internal logic of \alttext{\( \cE \)}{E}}
\subsection{Internal Kuroda = external Borel} \label{subsec:internal-external}
We recall the Kuroda double-negation translation;
see \cite{ferreira-oliva-double-negation}
for an overview of double-negation translations.

For a formula \(\phi\),
define the formula \(\phi^\ku\)
by replacing every instance of \(\forall x \, \psi\)
by \(\forall x \, \lnot \lnot \psi\)
(this is a recursive definition).
The formula \(\lnot \lnot \phi^\ku\) is called the
\defword{Kuroda (double-negation) translation} of \(\phi\).
Kuroda showed (see \cite{kuroda-double-negation})
that a first-order formula is a classical validity
iff its Kuroda translation is an intuitionistic validity.

\begin{defn} \label{defn:toa}
  We define the language \( \cL_{\TOA} \) of \emph{third-order arithmetic} as follows:
  \begin{itemize}
  \item \( \cL_{\TOA} \) has five sorts \( \N \), \( \N^\N \), \( (\N^\N)^{(\N^\N)} \), \( (\N^\N)^\N \), and \( \N^{(\N^\N)} \).
  \item For each sort of the form \( B^A \), there is a binary function symbol \( \ev = \ev_{A,B} \) of arity \( \ev \colon B^A \times A \to B \).
  \item There is a constant symbol \( 0 \tolon \N \) and a unary function symbol \( s \) of arity \( s \colon \N \to \N \).
  \end{itemize}
\end{defn}

\begin{rmk}
  Note that we have omitted power sorts from \( \cL_{\TOA} \).
  The reason is that we are only interested in the \emph{decidable} predicates, and the decidable predicates of sort \( \pow(X) \) are (provably in intuitionistic higher-order logic) in bijection with the sort \( \two^X \), which we have.
  If we were to include the sorts \( \pow(X) \), we would have to make a further classically invisible translation in Theorem~\ref{thm:internal-external} below, in addition to the Kuroda translation, in which we relativize the quantifiers over all power sorts to range over the decidable predicates only.

  We will discuss other extensions of \( \cL_{\TOA} \) in \S\ref{subsubsec:3oa-extensions} below.
\end{rmk}

Besides the interpretation of \( \cL_{\TOA} \) in any elementary topos \( \cD \) with natural numbers object, using exponential objects, there is in addition a natural \emph{Borel} interpretation, which differs from the obvious set-theoretic interpretation (i.e., the interpretation in the topos \( \cD = \Set \)) only in that, with \( A \) either of \( \N \) or \( \N^\N \), we interpret \( A^{(\N^\N)} \) as the set of \emph{Borel} maps \( \N^\N \to A \).

As above, let \( \cD \) denote either of \( \cE \) or \( \cE^\inc \).
We will show the following:
\begin{thm}[Internal Kuroda = external Borel]\label{thm:internal-external}
  Let \( \phi \) be a sentence in \( \cL_{\TOA} \).
  The Borel interpretation of \( \phi \) holds iff \( \cD \Vdash \lnot \lnot \phi^\ku \).
\end{thm}
\begin{proof}
  This follows from the more general Theorem~\ref{thm:internal-external-general} below.
\end{proof}
\noindent
By the soundness of the internal logic of a topos with respect to intuitionistic higher-order arithmetic (see \S\ref{subsec:elementary-topoi}), we conclude:
\begin{cor} \label{cor:hoa-implies-borel}
  Given a sentence \( \phi \) in \( \cL_{\TOA} \), if \( \neg \neg \phi^\ku \) is a theorem of intuitionistic higher-order arithmetic, then the Borel interpretation of \( \phi \) holds.
   \comm{Note: if \( \phi \) is classically valid, then \( \neg \neg \phi^\ku \) is valid; it seems that the same is true with ``valid'' replaced by ``theorem of higher-order arithmetic'', but this doesn't buy us anything here.
       The converse would buy us something, but seems to be false: indeed, there are surely things which are Borel true but no classically true.}
\end{cor}
\noindent
We give several examples in \S\ref{subsec:borel-examples} below which serve to illustrate the general principle that statements in Borel mathematics tend to be expressible in \( \cL_\TOA \), and hence subject to the above Theorem~\ref{thm:internal-external}.

We will deduce Theorem~\ref{thm:internal-external} from the following more general statement.
Given a small subcanonical presite \( (\cC,J) \), let \( \cL_{\cC} \) be the language with sorts \( B^A \) for each \( A,B \in \cC \), with function symbols \( \circ = \circ_{A,B,C} \colon C^B \times B^A \to C^A \) for \( A,B,C \in \cC \), and with a constant \( \bar f \tolon B^A \) for each \( f \colon A \to B \) in \( \cC \).
Again, we have two natural interpretations of \( \cL_{\cC} \): an interpretation in the internal language \( \cL_{\Sh(\cC, J)} \) in which \( B^A \) is interpreted as the exponential object \( \wh B^{\wh A} \) between representables; and an ordinary set-theoretic interpretation in which \( B^A \) is interpreted as the set \( \cC(A,B) \) of morphisms \( A \to B \) in \( \cC \).

We say that an object \( X \in \cC \) is \defword{empty} if \( X \) is covered by the empty family; or equivalently, by Lemma~\ref{lem:empty-object}, if \( \wh X \in \Sh(\cC, J) \) is initial.
Note that, by subcanonicty, the terminal object \( \tm \in \cC \) is non-empty.
We now have:
\begin{thm}\label{thm:internal-external-general}
  Let \( (\cC,J) \) be a small presite, and suppose that \( J \) is subcanonical and that each non-empty \( X \in \cC \) admits at least one morphism \( \tm \to X \) from the terminal object.

  Then a sentence \( \phi \) in \( \cL_{\cC} \) holds with respect to morphisms in \( \cC \) if and only if \( \Sh(\cC, J) \Vdash \lnot \lnot \phi^\ku \).
\end{thm}

We note that Theorem~\ref{thm:internal-external} follows immediately from Theorem~\ref{thm:internal-external-general}, since the language \( \cL_{\TOA} \) embeds into \( \cL_{\SBor} \) and---by Corollary~\ref{cor:nat-is-nat} and Proposition~\ref{propn:exps-are-exps}---the two interpretations considered in Theorem~\ref{thm:internal-external-general} restrict to the two interpretations considered in Theorem~\ref{thm:internal-external}.

\begin{lem}
  Let \( (\cC,J) \) be as in Theorem~\ref{thm:internal-external-general}, and let \( \phi \) be a formula in \( \cL_{\Sh(\cC,J)} \) with free variables among \( \vx \tolon \vX \).
  \begin{enumerate}
  \item
    Given any \( \tm \to \prod_i X_i \) and non-empty \( U \in \cC \),
    we have \( \tm \Vdash \phi \) iff \( U \Vdash \phi \).
  \item For any \( \wh U \to \prod_i X_i \), we have \( U \Vdash \neg \neg \phi \) iff \( \tm \Vdash \phi \) for all \( \tm \to U \).
  \end{enumerate}
\end{lem}
\begin{proof}
  The first claim is immediate from the assumption that any non-empty \( U \) admits a morphism \( \tm \to U \).

  For the second claim, assuming \( U \Vdash \neg \neg \phi \) and given \( \tm \to U \), we have \( \tm \Vdash \neg \neg \phi \), meaning that for every non-empty \( U' \to \tm \), there is a non-empty \( U'' \to U' \) with \( U'' \Vdash \phi \).
  In particular, taking \( U' = \tm \), we conclude that \( U'' \Vdash \phi \) for some non-empty \( U'' \), and hence \( \tm \Vdash \phi \) by part the first part of the lemma.

  Conversely, assuming \( \tm \Vdash \phi \) for all \( \tm \to U \), and given \( U' \to U \) non-empty, we must find \( U'' \to U' \) non-empty with \( U'' \Vdash \phi \).
  But by assumption, there exists a morphism \( \tm \to U' \), so we can just take \( U'' = \tm \).
\end{proof}

Theorem~\ref{thm:internal-external-general} follows immediately from the following inductive description of \( \tm \Vdash \phi \) in \( \Sh(\cC,J) \):
\begin{propn}\label{propn:main-thm-main-lemma}
  Let \( (\cC,J) \) be as in Theorem~\ref{thm:internal-external-general}, let \( \phi \) and \( \psi \) be formulas in \( \cL_{\Sh(\cC,J)} \) with free variables among \( \vx \tolon \vX \), and fix a morphism \( \tm \to \prod_i X_i \).
  The following hold in \( \Sh(\cC,J) \):
  \begin{enumerate}
    \item
      \(\tm \Vdash \top\).
    \item
      \(\tm \not \Vdash \bot\).
    \item
      \(\tm \Vdash \phi \land \psi\)
      iff \(\tm \Vdash \phi\) and \(\tm \Vdash \psi\).
    \item
      \(\tm \Vdash \bigvee_i \phi_i\)
      iff there is some \(i\) with \(\tm \Vdash \phi_i\).
    \item
      \(\tm \Vdash \phi \to \psi\)
      iff \(\tm \Vdash \phi\) implies \(\tm \Vdash \psi\).
    \item
      \(\tm \Vdash \exists x \tolon X.\ \phi\)
      iff there is some \(\tm \to X\)
      such that \(\tm \Vdash \phi\).
    \item
      \(\tm \Vdash \forall x \tolon X.\ \lnot \lnot \phi\) iff
      for all \(\tm \to X\),
      we have \(\tm \Vdash \phi\).
  \end{enumerate}
\end{propn}

\begin{proof}
  \leavevmode
  \begin{enumerate}
    \item[(1)-(3)]
      Immediate from the Proposition~\ref{propn:kripke-joyal}~\ref{item:forcing-conditions}.
      \setcounter{enumi}{3}
    \item
      We have \(\tm \Vdash \bigvee_i \phi_i\)
      iff there is a cover \( \bA \in J(\tm) \)
      such that for every \( a \in \bA \),
      we have \( U_a \Vdash \phi_i \) for some \( i \).
      In one direction,
      if \(\tm \Vdash \phi_i\),
      then we may take \( \bA \) to be the trivial cover \( \set{\id_\tm} \).
      In the other direction,
      given \( \bA \in J(\tm) \), some \( U_a \) must be non-empty since \( \tm \) is non-empty,
      so since \(U_a \Vdash \phi_i\) for some \( i \),
      we have \(\tm \Vdash \phi_i\) by the lemma.
    \item
      We have \(\tm \Vdash \phi \to \psi\)
      iff for all \(U \to \tm\),
      we have \(U \Vdash \phi\) implies \(U \Vdash \psi\).
      When \(U\) is empty,
      we have \(U \Vdash \psi\),
      and when \(U\) is non-empty,
      this is equivalent to \(\tm \Vdash \phi\) implies \(\tm \Vdash \psi\) by the lemma,
      so we are done.
    \item
      We have \(\tm \Vdash \exists x \, \phi\)
      iff there is some \( \bA \in J(\tm) \) and morphisms \( U_a \to X \) for \( a \in \bA \)
      such that \( U_a \Vdash \phi \) for all \( a \).
      If this holds, then some \( U_a \) must be non-empty, and thus \(\tm \Vdash \phi\) by the lemma.
      Conversely,
      if there is some \(\tm \to X\) so that \(\tm \Vdash \phi\),
      then we can just take \( \bA \) to be the trivial cover \( \set{\id_\tm} \).
    \item
      We have \(\tm \Vdash \forall x \, \lnot \lnot \phi\)
      iff for all \(U \to X\),
      we have \(U \Vdash \lnot \lnot \phi\),
      iff (by the lemma) for all \(\tm \to U\) and \(U \to X\),
      we have \(\tm \Vdash \phi\),
      iff for all \(\tm \to X\),
      we have \(\tm \Vdash \phi\).
      \qedhere
  \end{enumerate}
\end{proof}

\begin{proof}[\textbf{\textup{Proof of Theorem~\ref{thm:internal-external-general}}}]
  The claim follows by induction using Proposition~\ref{propn:main-thm-main-lemma}.
  The base case involves a further easy induction, showing that for a term \( t \tolon B^A \) in \( \cL_{\cC} \) with free variables among \( \vx \) where \( x_i \tolon {B_i}^{A_i} \), and given \( f_i \in \cC(A_i, B_i) \cong \Sh(\cC, J)(\tm, \wh{B_i}^{\wh A_i}) \), the interpretation \( \sem{t}_{\vx,\vec{f}} \colon \tm \to \wh B^{\wh A} \) of \( t \) in \( \Sh(\cC, J) \) agrees with the one in terms of morphisms in \( \cC \) (with respect to the bijection \( \cC(A, B) \cong \Sh(\cC, J)(\tm, \wh B^{\wh A}) \)).
\end{proof}

\comm{
  We had meant to put a remark here pointing out that Theorem~\ref{thm:internal-external-general} in fact holds for the \emph{presheaf} category \( \wh{\SBor} \), but that Theorem~\ref{thm:internal-external} does not follow in this case, since the interpretation of \( \cL_{\TOA} \) into \( \wh \SBor \) coming from the former is not the interpretation appearing in the latter (since \( \wh \N \in \wh{\SBor} \) is no longer the NNO in this case).

  However, now that Theorem~\ref{thm:internal-external-general} refers to \emph{non-empty}, rather than \emph{initial} objects (as I think it must), it seems to me that it no longer applies to the presheaf category: in the topology on \( \SBor \) giving the presheaf category (namely, the trivial topology, in which only split epis are covers), the initial object \emph{is} non-empty, but doesn't admit a map from \( \tm \).
}

\subsubsection{Definitional extensions} \label{subsubsec:3oa-extensions}
Theorem~\ref{thm:internal-external} establishes the agreement of the Borel interpretation and the interpretation in the Borel topos for sentences in the minimal third-order arithmetic language of Definition~\ref{defn:toa}.
This suffices in principle to express all of the Borel mathematics statements of interest, however in practice it is much more convenient to allow oneself certain \emph{definitional extensions} of this language (see \S\ref{subsec:elementary-topoi}).

The extensions which are available here are obtained as follows.
To begin with, we have \emph{second-order arithmetic} (2OA), with sorts \( \N \) and \( \N^\N \), both of which can be interpreted in the category \( \SBor \).
We then consider an arbitrary definitional extension \( \cL \) of 2OA which is \emph{still supported by} (i.e., can be interpreted in) \( \SBor \).
We then finally extend \( \cL \) by introducing exponential sorts between any two of its sorts.
This can no longer be interpreted in \( \SBor \) directly, but can still allows the interpretation using the hom sets of \( \SBor \).

We will describe a particular such extension which is convenient, suffices for all our purposes, and serves to illustrate the general principle, but is otherwise quite arbitrary.
We briefly comment on a more comprehensive class of extensions in Remark~\ref{rmk:infinitary-extensions} below.

\begin{defn} \label{defn:soa-toa-plus}
  We define the languages \( \cL_{\SOA}^+ \) and \( \cL_{\TOA}^+ \) as follows.
  The language \( \cL_\SOA^+ \) has
  \begin{itemize}
  \item the base sorts \( \N \) and \( \N_{<n} \) for each natural number \( n \ge 0 \), and
  \item sorts \( A \times B \), \( A + B \), \( B^A \), and \( A^{<\N} \) for any sorts \( A \) and \( B \),
  \end{itemize}
  and the following function symbols:
  \begin{itemize}
  \item the constant \( 0,1 \) of sort \( \N \) and the function symbols \( +,\times \colon \N \times \N \to \N \),
  \item the constant symbols \( \und 0, \ldots, \und {n-1} \) of sort \( \N_{<n} \) for each \( n \),
  \item the function symbols \( A \xot{\pi_A} A \times B \tox{\pi_B} B \) and \( A \tox{i_A} A + B \xot{i_B} B \) for each pair of sorts \( A \) and \( B \),
  \item the function symbol \( \br{-,-} \) of arity \( A \times B \to A \times B \) and \( \ev_{A,B} \) of arity \( B^A \times A \to B \) for each \( A \) and \( B \),
  \item the constant \( \nil_A \colon A^{<\N} \) and function symbol \( \cons_A \colon A \times A^{<\N} \to A^{<\N} \) for any sort~\( A \).
  \end{itemize}
  We then define \( \cL_{\TOA}^+ \) to be the extension of \( \cL_\SOA^+ \) obtained by introducing a sort \( B^A \) and a function symbol \( \ev_{A,B} \) as above for each pair of sorts \( A \) and \( B \) in \( \cL_\SOA^+ \).
\end{defn}
\noindent
The language \( \cL_{\TOA}^+ \) is simply a fragment of a definitional extension of higher-order arithmetic, in the sense of \S\ref{subsec:elementary-topoi}, and thus admits an interpretation into any elementary topos (for more on the interpretation of the sorts \( A^{<\N} \), which are interpreted as \emph{list objects}, see \S\ref{subsec:borel-myopia}).
And again, we have a second natural interpretation of \( \cL_\TOA^+ \), in which each sort in \( \cL_{\SOA}^+ \) is interpreted as the set of points in a standard Borel space, and hence each sort in \( \cL_{\TOA}^+ \) as the set of Borel maps between two standard Borel spaces.

The proofs of Theorem~\ref{thm:internal-external} and Corollary~\ref{cor:hoa-implies-borel} then carry over to give:
\begin{thm} \label{thm:internal-external-plus}
  Let \( \phi \) be a sentence in \( \cL_{\TOA}^+ \).
  The Borel interpretation of \( \phi \) holds iff \( \cD \Vdash \lnot \lnot \phi^\ku \).
\end{thm}
\begin{cor}
  Given a sentence \( \phi \) in \( \cL_{\TOA}^+ \), if \( \neg \neg \phi^\ku \) is a theorem of intuitionistic higher-order arithmetic, then the Borel interpretation of \( \phi \) holds.
\end{cor}

\begin{proof}[\textbf{\textup{Proof of Theorem~\ref{thm:internal-external-plus}}}]
  The proof is the same as that of Theorem~\ref{thm:internal-external}, except that Theorem~\ref{thm:internal-external-general} doesn't quite apply directly, since there isn't an interpretation of \( \cL_{\SOA}^+ \) into \( \cL_{\SBor} \).
  However, this can be corrected via a small modification to the language \( \cL_{\cC} \) appearing in Theorem~\ref{thm:internal-external-general}: in addition to the sorts \( B^A \) and the function symbols \( C^B \times B^A \to C^A \), we also include function symbols \( {A_1}^\tm \times \cdots \times {A_n}^\tm \to B^\tm \), with \( \tm \) the terminal object, for each morphism \( f \colon \prod_{i=1}^n A_i \to B \) in \( \cC \), as in the usual internal language of a category (see \S\ref{subsec:elementary-topoi}).
  The two interpretations of \( \cL_{\cC} \) appearing in Theorem~\ref{thm:internal-external-general} extend to this larger language, and the same proof, using Proposition~\ref{propn:main-thm-main-lemma}, shows that theorem still holds with this modified definition of \( \cL_{\cC} \).
  This then immediately implies our claim, using the interpretation of \( \cL_{\SOA}^+ \) into \( \cL_{\SBor} \) which is now available, given the modifications made to \( \cL_{\SBor} \).
\end{proof}

\begin{rmk} \label{rmk:infinitary-extensions}
  A more systematic description of the available definitional extensions of \( \cL_{\TOA} \) in the present context would involve the use of (countably) infinitary logic, so that, rather than just introducing exponential sorts \( X^\N \), we could more generally introduce countable products \( \prod_n X_n \) (and coproducts \( \coprod_n X_n \)) of sorts, and the sort \( X^{<\N} \) could be defined as a sub-sort of  \( X^\N \), given by the formula \( \phi(x) \equiv \bigvee_{k \in \N} \bigwedge_{n \in \N} \pi_n(x) = \pi_k(x) \).
\end{rmk}

\subsection{Examples from Borel mathematics} \label{subsec:borel-examples}
We now examine some well-known statements in Borel mathematics, observe that they can be formulated in third-order arithmetic (that is, in \( \cL_{\TOA}^+ \)---see \S\ref{subsubsec:3oa-extensions}), and conclude that they hold in the Borel topoi \( \cE \) and \( \cE^\inc \).

\subsubsection{Continuum hypothesis}
The continuum hypothesis is classically equivalent
to various inequivalent intuitionistic statements,
some of which hold in \(\cD\) and some of which do not.

For instance,
consider the following version:
``for all \(f : \N^\N \to 2\),
either there exists \(g : \N \to \N^\N\) such that
for all \(x \in \N^\N\) satisfying \(f(x) = 0\)
there exists \(n \in \N\) satisfying \(g(n) = x\),
or there exists \(h : \two^\N \hto \N^\N\) such that
for all \(\alpha \in \two^\N\),
we have \((f \circ h)(\alpha) = 0\).''

This doesn't hold in \(\cD\),
since this would be saying that for every standard Borel space \(U\)
and every Borel subset \(R\) of \(U \times \N^\N\),
the set \(\{u \in U : \text{\(R_u\) is countable}\}\) is Borel,
which is not always the case.
But since this holds for Borel maps by the perfect set theorem \cite[Theorem~29.1]{kechris-classical},
its Kuroda translation holds in \(\cD\).

\subsubsection{Borel determinacy}
Consider the sentence \(\phi_{\text{det}}\) which says the following:
for every \(f : \N^\N \to 2\),
there exists \(i : 2\) and \(\sigma : \N^{<\N} \to \N\)
such that for every \(\alpha : \N^\N\),
if every \(n : \N\) with \(n \equiv i \pmod 2\)
satisfies \(\alpha\uhr (n+1) = (\alpha\uhr n) \concat \sigma(\alpha\uhr n)\),
then \(f(\alpha) = i\).

Then \(\phi_{\text{det}}\) holds for Borel maps;
this is the Borel determinacy theorem of Martin
(see \cite[Theorem~20.5]{kechris-classical}).
Thus internally,
the Kuroda translation of \(\phi_{\text{det}}\)
holds in $\cD$.

However,
\(\phi_{\text{det}}\) does not hold in $\cD$,
since given a parametrized family of games,
the set where Player I wins is not necessarily Borel
(it can be any analytic set).

\subsubsection{Well-foundedness of the Wadge hierarchy}
Given two subsets $A$ and $B$ of $\N^\N$,
a \textbf{Wadge reduction} from $A$ to $B$
is a continuous function $\N^\N \to \N^\N$
under which the preimage of $B$ is $A$,
and we say that $A$ is \textbf{Wadge reducible} to $B$,
denoted $A \le_W B$,
if there is a Wadge reduction from $A$ to $B$.
Wadge reducibility is a preorder on subsets of $\N^\N$,
and the Martin-Monk theorem says that its restriction to Borel subsets is well-founded
\cite[Theorem~21.15]{kechris-classical}.
Classically,
assuming Dependent Choice,
this is equivalent to saying that for all sequences $(A_n)_n$ of Borel subsets of $\N^\N$,
there is some $n \in \N$ such that $A_{n+1}$ does not Wadge reduce to $A_n$.
In other words,
the following sentence's Kuroda translation holds in the Borel topos:
``for all functions $f : \N \times \N^\N \to 2$,
there is some $n \in \N$ such that for all continuous functions $g : \N^\N \to \N^\N$,
we have $f_{n+1} \neq f_n \circ g$.''
Note that continuity is expressible in this language:
a function $\N^\N \to \N^\N$ is continuous
iff its compositions with the projections are continuous,
and a function $\N^\N \to \N$ is continuous
iff every element of $\N^\N$ has an initial segment
such that every element of $\N^\N$ with that initial segment maps to the same value.

\subsubsection{The $E_0$ dichotomy}
Consider the equivalence relation $E_0$ on $\two^\N$ defined as follows:
$x$ and $y$ are $E_0$-related if there exists a natural number $n$
such that for all natural numbers $k > n$,
we have $x_k = y_k$.

Consider the statement
``for every equivalence relation $E$ on $\N^\N$,
either there is an injection from $\N^\N/E$ to $\two^\N$,
or there is an injection from $\two^\N/E_0$ to $\N^\N/E$.''
The Kuroda translation of this holds in the Borel topos:
the external version of this statement is the $E_0$ dichotomy
of Harringon-Kechris-Louveau
\cite{HKL90}.

Similarly,
    we can also consider Silver's theorem \cite[Theorem~35.20]{kechris-classical},
which says that every Borel equivalence relation
either has countably many classes
or a perfect set of classes.

\subsubsection{Borel combinatorics}
Consider the statement
``every $d$-regular graph $(\N^\N)^2 \to 2$ on $\N^\N$ has a proper $(d+1)$-coloring''.
We are viewing a function $(\N^\N)^2 \to 2$ as a graph with vertex set $\N^\N$,
where the set of edges is the preimage of $1$,
a graph is $d$-regular if every vertex has $d$ neighbors,
and a proper $(d+1)$-coloring is a function $\N^\N \to d+1$
such that adjacent vertices receive different colors.
Its Kuroda translation holds in the Borel topos,
since externally,
every $d$-regular Borel graph on $\N^\N$ has a Borel proper $(d+1)$-coloring
(\cite{KST99}).

On the other hand,
consider the statement
``every $d$-regular acyclic graph $(\N^\N)^2 \to 2$ on $\N^\N$ has a proper $d$-coloring'',
where a graph is acyclic if it has no cycles.
Its Kuroda translation does not hold in the Borel topos,
since the corresponding external statement is that
every $d$-regular acyclic Borel graph on $\N^\N$ has a Borel proper $d$-coloring,
which is false \cite{Mar16}.

\subsubsection{Choice for families of countable sets} \label{subsubsec:countable-choice}
Consider the statement
``for every surjection \( f \colon \two^\N \to \two^\N \),
  if for every \( y \in \two^\N \)
  there is \( h \colon \two^\N \to \N \)
  which is injective on the \( f \)-fibre over \( y \),
  then there is \( g \colon \two^\N \to \two^\N \) with \( f \circ g = \id_\N \).''
This is a form of the axiom of choice for families of countable sets.
Its Borel interpretation holds by the Lusin-Novikov theorem
\cite[Theorem~18.10]{kechris-classical},
and hence its Kuroda translation holds in the Borel topos.
Actually,
even without Kuroda translation,
it already holds in the Borel topos.
To see this,
fix \( U_0 \) and a Borel \( f \colon U_0 \times \two^\N \to \two^\N \),
and assume that \( U_0 \) forces
``\( f \) is surjective
and for every \( y \in \two^\N \),
there is \( h \colon \two^\N \to \N \)
injective on its \( f \)-fibre''.
We need to show that \( U_0 \) forces
``there is a section \( g \) of \( f \).''

Now \( f \) induces a function \( U_0 \times \two^\N \to \two^\N \),
and we claim that every fibre of the induced function
\( U_0 \times \two^\N \to U_0 \times \two^\N \)
is countable.
To see this,
let \( z \in \two^\N \).
Since \( U_0 \) forces ``\( f \) is surjective'',
there is some cover \( U_1 \to U_0 \)
and some Borel \(x \colon U_1 \to \two^\N \)
such that \( U_1 \) forces ``\( f(x) = z \)''.
Now since \( U_0 \) forces
``for every \( y \in \two^\N \), there is \( h \colon \two^\N \to \N \)
injective on its \( f \)-fibre'',
there is some cover \( U_2 \to U_1 \)
and some Borel \( h \colon U_2 \times \two^\N \to \N \)
such that \( U_2 \) forces
``\( h \) is injective on the \( f \)-fibre of \( z \)''.
Now let \( u \in U_2 \).
We claim that the fibre of \( (u, z) \) under the map
\( U_2 \times \two^\N \to U_2 \times \two^\N \) induced by \( f \)
is countable,
which implies the result for \( U_0 \).
Suppose that \( (u, x_0) \) and \( (u, x_1) \) both map to \( (u, z) \),
and suppose that \( h(u, x_0) = h(u, x_1) \).
Consider \( x \colon U_2 \to \two^\N \),
and let \( x_0', x_1' \colon U_2 \to \two^\N \)
be the same function except sending
\( u \) to \( x_0 \) and \( x_1 \) respectively.
Then \( U_2 \) forces ``\( f(x_0') = f(x_1') = z \)''
and ``\( h(x_0') = h(x_1') \)'',
so since \( U_2 \) forces ``\( h \) is injective on the fibre'',
it forces ``\(x_0' = x_1'\)'',
which in particular means that \( (u, x_0) = (u, x_1) \),
as desired.

So by the (external) Lusin-Novikov theorem,
this has a Borel section,
and this witnesses that \( U' \) forces that
``there is a section \( g \) of \( f \)''.

By contrast, the axiom of choice in general fails in the Borel topos; see \S\ref{subsec:more-topos-props}.

\subsection{Borel myopia} \label{subsec:borel-myopia}
The Church-Turing thesis is the informal statement that the computable functions on the natural numbers are exactly those whose values can be computed in a finite time by a deterministic procedure.
A striking application of elementary topos theory, and more generally intuitionistic set theory, is that this informal statement can be given a formal counterpart: there are topoi in which the statement ``every function \( \N \to \N \) is computable'' (which is obviously false in classical set theory) holds internally. \cite{hyland-effective-topos,nashaat-powell-valued-models}

In the present context of a Borel topos, it is natural to ask whether the corresponding question, ``is every function \( \R \to \R \) Borel?'' holds internally; we will express this by saying the topos is Borel-myopic.
Clasically, this question is equivalent to ``is every subset \( X \subset \R \) Borel?''.
Intuitionistically, the situation is more subtle.
First of all, rather than \( \R \), we will ask the question about \( \two^\N \) (which may no longer be equivalent).
Then, the statement that every function \( \two^\N \to \two^\N \) is Borel turns out to be equivalent to every \emph{decidable} subset of \( \two^\N \) being Borel.
Indeed, the forward implication is clear, since we can take an indicator function \( \two^\N \to \two \hto \two^\N \) on any decidable subset.
And the reverse implication holds because a map \( \two^\N \to \two^\N \) is Borel as soon as each of its components \( \two^\N \to \two \) is.\footnote{
    There is a similar circumstance in the setting of computability: that every function \( \N \to \N \) is computable is equivalent, intuitionistically, to every \emph{decidable} subset of \( \N \) being computable.
}
We will take these equivalent statements---that all functions \( \two^\N \to \two^\N \), or all decidable subsets \( A \subseteq \two^\N \), are Borel---to be our definition of Borel myopia.

We must now address a potential ambiguity in what we might mean by ``Borel set''---an ambiguity which does not arise classically, but does in the intuitionistic context, and more generally in the absence of the axiom of choice.
The most immediate definition of Borel set is: member of the least \( \sigma \)-algebra of \( \pow(\two^\N) \) containing the \emph{basic opens}---i.e., the sets of sequences with a fixed initial segment.
However, we can also consider subsets of \( \two^\N \) admitting a \emph{Borel code}:
a countably branching labeled tree encoding the way in which the given set is built out of basic opens using complements and countable unions and intersections.
Using the axiom of choice, one can always reconstruct a Borel code from any given Borel set, by presenting it as a negation or countable union or intersection, and recursively choosing a Borel code for each of the constituent Borel sets.
But in the absence of choice, the two concepts are distinct.
We call a Borel set admitting a Borel code \defword{uniformly Borel}; obviously (even without choice, or excluded middle), any uniformly Borel set is Borel.

The significance of this distinction for us is that, if we take the weaker notion of Borel, there are models of set theory in which \emph{every} subset---not just every decidable subset---of \( \two^\N \) is Borel.
We might call this property \defword{strong Borel myopia}.
For example, there exists a model of ZF due to Feferman and Levy in which \( \two^\N \) is a countable union of countable sets \( A_n \) \cite[Theorem~10.6]{jech-axiom-of-choice}.
Hence, for any \( X \subset \two^\N \), we have that \( X = \bigcup_{n \in \N} (X \cap A_n) \) is itself a countable union of countable sets, and hence Borel.
However, \( X \) needn't be \emph{uniformly} Borel, since this would require enumerating each of the sets \( X \cap A_n \).\footnote{
  Presumably, there are not only models of ZF, but also topoi in which \( \two^\N \) is a countable union of countable sets, and hence in which every subset is Borel, but we have not checked this.
}

Indeed, as we show in Remark~\ref{rmk:no-strong-myopia} below, there is no topos satisfying \emph{strong uniform} Borel myopia, i.e., in which \emph{every} subset of \( \two^\N \) is uniformly Borel.
The reason is simply that we can (intuitionistically) explicitly construct a set which is not uniformly Borel.
We thus have the following four version of Borel myopia, which are related by the displayed implications, and which have the displayed models:
\begin{center}
  \underline{Types of Borel myopia}\\[0.5em]
  \begin{tabular}{c|c|c}
    & Plain & \( \Leftarrow \) Uniform
    \\
    \hline
    Decidable & Feferman-Levy, \( \cE \) & \( \cE \) \\
    \hline
    \shortstack{%
    \( \Uparrow \) \\
    Strong
    } & Feferman-Levy & None
  \end{tabular}
\end{center}
There are two further properties which we want to consider, namely the \emph{converses} of decidable Borel myopia: that every Borel set, resp., every uniform Borel set, is decidable.
It turns out that these both hold in \( \cE \):
\begin{center}
  \underline{Converse to Borel myopia}\\[0.5em]
  \begin{tabular}{c|c}
    All plain Borel sets are decidable \( \Rightarrow \) & All uniform Borel sets are decidable
    \\
    \hline
    \( \cE \) & \( \cE \)
  \end{tabular}
\end{center}
(Of course, the converse also holds in Feferman-Levy, for the trivial reason that it satisfies excluded middle.)
Thus, to summarize, in \( \cE \), the Borel, uniform Borel, and decidable subsets of \( \two^\N \) all coincide.

We now turn to the expression of these properties in a topos---which is to say, in higher-order arithmetic.
For this purpose, it will be useful to work with a certain definitional extension of higher-order logic (see \S\ref{subsec:elementary-topoi}).
Specifically, we consider the language with base sorts \( \N \) and \( \N_{<\N} \) and with sort-formers \( A \times B \), \( A + B \), \( B^A \), and \( A^{<\N} \) as in Definition~\ref{defn:soa-toa-plus}, as well as \( \pow(A) \).
The function and relation symbols are all the ones appearing in Definition~\ref{defn:soa-toa-plus} as well as the element relations \( \in_A \) of arity \( A \times \pow(A) \) for each sort \( A \).
We leave it to the reader to give the relevant formulas defining each of these sorts and relation and function symbols as in \S\ref{subsec:elementary-topoi}.
We will introduce a few more definable functions into our extension below.

Being a definitional extension of higher-order arithmetic, this language admits an interpretation into any elementary topos with NNO \( \cC \).
Specifically, the various sorts are interpreted as using the NNO, products, coproducts (with \( \N_{<\N} \) the \( n \)-fold coproduct of the terminal object), exponentials, and power objects of \( \cC \), and using \emph{list objects} to interpret the sorts \( A^{<\N} \).

Given an object \( X \) in a category \( \cC \), a \defword{list object} on \( X \) is an object \( X^{<\N} \) equipped with morphisms \( \tm \tox{\nil} X^{<\N} \xot{\cons} X \times X^{<\N} \) and satisfying the evident universal property: given any other object \( Y \) together with morphisms \( \tm \tox{a} Y \xot{f} X \times Y \), there is a unique morphisms \( g \) making the following diagram commute.
  \[
    \begin{tikzcd}
      \tm \ar[d, equals] \ar[r, "\nil"] & X^{<\N} \ar[d, "g", dashed] & X \times X^{<\N} \ar[d, "\id \times g"] \ar[l, "\cons"'] \\
      \tm \ar[r, "a"] & Y & X \times Y. \ar[l, "f"']
    \end{tikzcd}
  \]
  In an elementary topos with natural numbers object, every object admits a list object \cite[Theorem~A2.5.17]{johnstone-elephant-vols}.

It will be convenient to further give ourselves function symbols \( \bigcap, \bigcup \colon \pow(X)^A \to \pow(X) \) and \( \neg \colon \pow(X) \to \pow(X) \); these are obviously definable as well, namely by
\[
  \begin{split}
    \bigcap f = p & \tot \pBig{\forall x \tolon X.\ x \in p \tot \pbig{\forall a \tolon A.\ x \in \ev(f, a)}} \\
    \bigcup f = p & \tot \pBig{\forall x \tolon X.\ x \in p \tot \pbig{\exists a \tolon A.\ x \in \ev(f, a)}} \\
    \neg q = p & \tot \pBig{\forall x \tolon X.\ x \in p \tot \neg x \in q}.
  \end{split}
\]
And finally, we introduce a function symbol \( i_{\basic} \colon \two^{<\N} \to \pow(\two^\N) \) corresponding to the definition
\[
  \begin{split}
    \forall x \tolon \two^\N.&\
    x \in i_{\basic}(\nil)
    \ \wedge\\
    & \forall n \tolon \two.\
    \forall a \tolon \two^{<\N}.\
    x \in i_{\basic} \pbig{\cons(n, a)} \tot \\
    & \qquad\pBig{
    \ev(x, 0) = n \wedge
    \exists b \tolon \two^\N.\
    b \in i_{\basic}(a) \wedge
    \forall m \tolon \N.\ \ev(b, m) = \ev\pbig{a, s(m)}
    }.
  \end{split}
\]
We will introduce a small further definitional extension below when we define uniform Borel myopia, but first, let us define plain Borel myopia:

\begin{defn}
  We define the formulas \( \cl_{\cap}(P), \cl_{\cup}(P) \) with a free variable of sort \( \pow\pbig{\pow(\two^\N)} \):
  \[
    \cl_\cap(P) \defequiv
    \forall Q \tolon \pow(\two^\N)^\N.\
    \pbig{\forall n \tolon \N.\
    \ev(Q, n) \in P} \to
    \bigcap Q \in P.
  \]
  The formula \( \cl_\cup(P) \) is defined the same way but with \( \bigcup \) in place of \( \bigcap \).
  Similarly, we define \( \cl_{\neg}(P) \) as
  \[
    \cl_{\neg}(P) \defequiv
    \forall q \tolon \pow(\two^\N).\
    q \in P \to \neg q \in p
  \]
  and \( \cl_{\basic}(P) \) as
  \[
    \cl_{\basic}(P) \defequiv
    \forall b \tolon \two^{<\N}.\
    i_{\basic}(b) \in P.
  \]
  We define the formula \( \Bo(p) \) with \( p \tolon \pow(\two^\N) \) to be
  \[
    \Bo(p) \defequiv
    \forall P \tolon \pow\pbig{\pow(\two^\N)}.\
    \pbig{
      \cl_{\cup}(P) \wedge \cl_{\cap}(P) \wedge \cl_{\neg}(P) \wedge \cl_{\basic}(P)
    } \to p \in P.
  \]

  We now say that a topos \( \cC \) satisfies \defword{plain Borel myopia} if
  \[
    \cC \vDash \forall p \tolon \pow(\two^\N).\
    (\forall x \tolon \two^\N.\ x \in p \vee x \notin p) \to
    \Bo(p)
  \]
  and that \( \cC \) satisfies the \defword{converse to plain Borel myopia} if
  \[
    \cC \vDash \forall p \tolon \pow(\two^\N).\
    \Bo(p) \to
    (\forall x \tolon \two^\N.\ x \in p \vee x \notin p).
  \]
\end{defn}

We now proceed to uniform Borel myopia, for which we will first of all need to define Borel codes.
As mentioned above, we will take a Borel code to be a labeled countably branching tree, which we interpret as a function from the set of finite paths from the root in \( \N^{<\N} \) to the set of labels.
A label is either a code for a basic open subset---thus, an element of \( \two^{<\N} \)---or a label indicating ``intersection'', ``union'', or ``complement''---thus, an element of \( \N_{<3} = \set{0,1,2} \).
Really, such a tree should be \emph{at most} countably branching: a node labeled with \( a \in \two^{<\N} \) should have no children, and a node labeled with ``\( 2 \)'' (``complement'') should have one child.
We adopt the convention that, in the first case, all descendants of the given node likewise carry the label \( a \), and in the second case, that all the children of the given node are equal.

Thus, it is not \emph{all} functions \( \two^{<\N} \to \N_{<3} + \two^{<\N} \) which serve as Borel codes, but only those satisfying the above two constraints.
However, we will not separately define the notion of ``Borel code'', but only that of ``coded Borel set''.

We define the sort
\[
  \code \defequiv
  (\N_{3} + \two^{<\N})^{\N^{<\N}}
\]
We now introduce four further function symbols: \( f_{\cap}, f_\cup \colon \code^\N \to \code \), \( \phi_\neg \colon \code \to \code \), and \( f_{\basic} \colon \two^{<\N} \to \code \), corresponding to the following definitions:
\[
  \begin{split}
    f_{\cap}(T) = t \tot\
    & \ev(t, \nil) = i_{\N_3}(0) \wedge
    \forall n \tolon \N.\ \forall a \tolon \N^{<\N}.\
    \ev\pbig{t, \cons(n, a)} = \ev\pbig{\ev(T, n), a}
    \\
    f_{\neg}(s) = t \tot\
    & \ev(t, \nil) = i_{\N_3}(2) \wedge
    \forall n \tolon \N.\ \forall a \tolon \N^{<\N}.\
      \ev\pbig{t, \cons(n, a)} = \ev(s, a)
    \\
    f_{\basic}(b) = t \tot\ &
    \forall a \tolon \N^{<\N}.\
    \ev(t, a) = i_{\two^{<\N}}(b).
  \end{split}
\]
The formula for \( f_{\cup} \) is the same as that for \( f_\cap \) but with \( \ev(t, \nil) = i_{\N_3}(1) \) in place of \( \ev(t, \nil) = i_{\N_3}(0) \).

\begin{defn}
  We define formulas \( \wt{\cl}_{\cap}(A),\wt{\cl}_{\cup}(A),\wt{\cl}_{\neg}(A),\wt{\cl}_{\basic} \) with \( A \) of sort \( \pow\pbig{\code \times \pow(\two^\N)} \):
  \[
    \begin{split}
      \wt{\cl}_{\cap}(A) & \defequiv
      \forall T \tolon \code^\N.\
      \forall P \tolon \pow(\two^\N)^\N.\
      \pbig{
      \forall n \tolon \N.\
      \br{T(n), P(n)} \in A
      } \to
      \br{f_{\cap}(T), \bigcap(P)} \in A
      \\
      \wt{\cl}_{\neg}(A) & \defequiv
      \forall t \tolon \code.\
      \forall p \tolon \pow(\two^\N)^\N.\
      (
      \br{t, p} \in A
      ) \to
      \br{f_{\neg}(t), \neg p} \in A
      \\
      \wt{\cl}_{\basic}(A) & \defequiv
       \forall b \tolon \basic.\
       \br{f_{\basic}(b), i_{\basic}(b)} \in A
    \end{split}
  \]
  (and with \( \wt{\cl}_{\cup} \) defined analogously to \( \wt{\cl}_{\cap} \)).
  We then define \( \cBo(q) \) (``\( q \tolon \code \times \pow(\two^\N) \) is a coded Borel set'') and \( \uBo(p) \) (``\( p \colon \pow(\two^\N) \) is a uniformly Borel set'') as follows:
  \[
    \begin{split}
      \cBo(q) & \defequiv
      \forall Q \colon \pow\pbig{\code \times \pow(\two^\N)}.\
      \pbig{
      \wt{\cl}_\cup(Q) \wedge
      \wt{\cl}_\cap(Q) \wedge
      \wt{\cl}_\neg(Q) \wedge
      \wt{\cl}_\basic(Q)
      } \to q \in A,
      \\
      \uBo(p) & \defequiv
      \exists t \colon \code.\ \cBo(\br{t, p}).
    \end{split}
  \]
  We then say that a topos \( \cC \) satisfies \defword{uniform Borel myopia} if
  \[
    \cC \vDash
    \forall p \tolon \pow(\two^\N).\
    (\forall x \tolon \two^\N.\ x \in p \vee x \notin p) \to
    \uBo(p)
  \]
  and satisfies \defword{the converse to uniform Borel myopia} if
  \[
    \cC \vDash
    \forall p \tolon \pow(\two^\N).\
    \uBo(p) \to
    (\forall x \tolon \two^\N.\ x \in p \vee x \notin p).
  \]
\end{defn}

Let us record the relation mentioned above between the plain and uniform variants:
\begin{propn} \label{propn:uniform-implies-plain}
  Any elementary topos with natural numbers object \( \cC \) satisfies
  \[
    \cC \vDash
    \forall p \tolon \pow(\two^\N).\
    \uBo(p)
    \to
    \Bo(p).
  \]
  Hence, if \( \cC \) satisfies uniform Borel myopia, then it satisfies plain Borel myopia, and if it satisfies the converse to plain Borel myopia, then it satisfies the converse to uniform Borel myopia.
\end{propn}

Now our main theorem is:
\begin{thm} \label{thm:uniform-borel-myopia}
  The Borel topos \( \cE \) satisfies uniform Borel myopia and the converse to plain Borel myopia; i.e.,
  \[
    \cC \vDash
    \forall p \tolon \pow(\two^\N).\
    \uBo(p)
    \tot
    \Bo(p)
    \tot
    \forall x \tolon \two^\N (x \in p \vee x \notin p).
  \]
\end{thm}

\begin{proof}[\textbf{\textup{Proof of Proposition~\ref{propn:uniform-implies-plain}}}]
  Let us show that the sentence in question is provable in intuitionistic higher-order arithmetic.
  By comprehension we have
  \[
    \vdash
    \exists Q \colon \pow\pbig{\code \times \pow(\two^\N)}.\
    \forall x \colon \code \times \pow(\two^\N).\
    x \in Q
    \tot
    \Bo\pbig{\pi_{\pow(\two^\N)}(x)}.
  \]
  Now, it is clear from the definitions of \( \Bo \) and the various \( \wt\cl \) formulas that
  \[
    \vdash
    \pBig{
      \forall x \colon \code \times \pow(\two^\N).\
      x \in Q
      \tot
      \Bo\pbig{\pi_{\pow(\two^\N)}(x)}
    } \to
    (\wt{\cl}_{\cup}(Q) \wedge \wt{\cl}_{\cap}(Q) \wedge \wt{\cl}_{\neg}(Q) \wedge \wt{\cl}_{\basic}(Q)).
  \]
  It thus follows from the definition of \( \cBo \) that
  \[
    \vdash
    \forall p \in \code \times \pow(\two^\N).\
    \cBo(p) \to \Bo\pbig{\pi_{\pow(\two^\N)}(p))},
  \]
  which implies
  \[
    \vdash
    \forall p \in \pow(\two^\N).\
    \pbig{\exists t \in \code.\ \cBo(\br{t, p})}
    \to
    \Bo(p),
  \]
  as desired.
\end{proof}

\begin{proof}[\textbf{\textup{Proof of Theorem~\ref{thm:uniform-borel-myopia}}}]
  To begin with, let us unfold, in the particular topos \( \cE \), all of the general definitions that went into the statement of Borel myopia:
  \begin{itemize}
  \item By Corollary~\ref{cor:nat-is-nat}, the natural numbers object in \( \cE \) is just the representable \( \wh \N \).
  \item By Proposition~\ref{propn:exps-are-exps}, the exponential ``\( \two^\N \)'' is just the representable \( \wh{\two^\N} \).
  \item The morphism \( i_{\basic} \colon \two^{<\N} \to \pow(\two^\N) \) corresponds to the representable subobject (i.e., Borel subset) of \( \two^{<\N} \times \two^\N \) consisting of \( (a, x) \in \two^{<\N} \times \two^\N \) such that \( a \) is an initial segment of \( x \).
    In particular, this is a \emph{decidable} subobject of \( \two^{<\N} \times \two^\N \) (by Proposition~\ref{propn:decidable-sub-is-borel}), hence \( i_{\basic} \) factors through \( \two^{\two^\N} \tto \pow(\two^\N) \).

    The morphism \( i_{\basic} \) is monic, and given \( U \in \SBor \), a morphism \( a \colon \wh U \to \pow(\two^\N) \) classifying a subobject \( A \tto \wh U \times \two^\N \) factors through the image \( i_\basic \) iff \( A \) is represented by a Borel subset \( B \subset U \times \two^\N \) of the following form: there is a countable disjoint cover \( U = \bigcup_n U_n \) and basic opens \( V_n \subset \two^\N \) so that \( B = \bigcup_n U_n \times V_n \).
  \item In \( \Set \), for any object \( X \), a list object \( X^{<\N} \) is given by the usual set of finite sequences in \( X \) (thus \( X^{<\N} \subset \pow(\N \times X) \)).
  Likewise, in \( \SBor \), list objects are given by the same construction, which is made into a Borel space by construing it as a countable coproduct \( \coprod_{n \ge 0} X^n \).
  Since products and countable coproducts are preserved by the Yoneda embedding \( \yon \colon \SBor \to \cE \), and since finite products distribute over countable coproducts in both \( \SBor \) and \( \cE \), it follows that \( \yon \) preserves list objects as well.

  Using that \( \yon \) also preserves exponentials by countable spaces (see Proposition~\ref{propn:sbor-lims-colims}), we conclude that the interpretation of the sort \( \code \) in \( \cE \) is representable, namely by the space \( \code = (\N_3 + \two^{<\N})^{\N^{<\N}} \in \SBor \); and the interpretation of \( \code^\N \) is likewise representable.
\item The morphisms \( f_{\basic} \colon \two^{<\N} \to \code \) and \( f_\neg \colon \code \to \code \) and \( f_\cap, f_\cup \colon \code^\N \to \N \) are given as follows:
  \( f_{\basic} \) takes \( a \in \two^{<\N} \) to the constant function \( x \mapsto (1, a) \in (\N_3 \sqcup \two^{<\N}) \); the function \( f_\neg \) takes \( t \in \code \) to the function taking \( \nil \in \N^{<\N} \) to \( (0, 2) \in (\N_3 \sqcup \two^{<\N}) \) and \( \cons(n, x) \in \N^{<\N} \) to \( t(x) \); the function \( f_\cap \) takes \( t \in \code^\N \) to the function taking \( \nil \in \N^{<\N} \) to \( (0, 0) \in (\N_3 \sqcup \two^{<\N}) \) and \( \cons(n, x) \in \N^{<\N} \) to \( t(n)(x) \); and \( f_\cup \) is defined similarly but with \( (0, 1) \) in place of \( (0, 0) \).
\item Fix \( U,X \in \SBor \) and a morphism \( r \colon \wh U \to \pow(\wh X)^\N \)---i.e., since \( \N \) is a countable coproduct of \( \tm \)'s in \( \cE \) (see Corollary\ref{cor:nat-is-nat}), and hence \( \pow(\wh X)^\N \) a countable product, a sequence of morphisms \( r_n \colon \wh U \to \pow(\wh X) \)---classifying subobjects \( A_n \tto \wh U \times \wh X \).
  Then \( \bigcap r, \bigcup r \colon \wh U \to \pow(\wh X) \) classify the subobjects \( A \tto \wh U \times \wh X \) given by the meet or join of the \( A_n \) in the subobject lattice of \( \wh U \times \wh X \).

  If each \( r_n \) is \emph{decidable}, i.e., factors through \( \two^\N \tto \pow(\wh X) \), so that \( A_n \to \wh U \times \wh X \) represents a Borel subset of \( U \times X \), then \( A \) is simply the intersection or union of these Borel subsets.
    This follows from Proposition~\ref{propn:subsheaves-in-e} and the fact that the embedding \( \Sub_{\SBor}(X) \to \AI(X) \) of Borel subsets of \( X \) into \( \sigma \)-ideals preserves countable meets and joins.
  \item Similarly, if \( r \colon \wh U \to \pow(\wh X) \) classifies a subobject (resp.\ decidable subobject) of \( \wh U \times \wh X \), then \( \neg r \colon U \to \pow(\wh X) \) is classified its complement in \( \Sub(\wh U \times \wh X) \) (resp.\ its complementary Borel subset of \( U \times X \)).
  \item For a \emph{global} element \( P \colon \tm \to \pow\pbig{\pow(\two^\N)} \) (the only case we will need) classifying a subobject \( A \tto \pow(\two^\N) \), we have that \( \tm \Vdash \cl_{\cap}(P) \) (in which case we say that \( A \) is \emph{closed under} \( \bigcap \)) if and only if, for each \( U \in \SBor \), and \( Q = (q_n)_{n \in \N} \colon \wh U \to \pow(\two^\N)^\N \), if \( q_n \circ u \colon U' \to \pow(\two^\N) \) factors through \( A \) for each \( n \in \N \) and each \( u \colon U' \to U \) in \( \SBor \), then \( \bigcap Q \colon U \to \pow(\two^\N) \) factors through \( A \) as well.
    The description of \( \cl_{\cup} \) is similar.
    And similarly, we have that \( A \tto \pow(\two^\N) \) is closed under \( \neg \) if and only if \( \neg \circ q \colon \wh U \to \pow(\two^\N) \) factors through \( A \) whenever \( q \colon \wh U \to \pow(\two^\N) \) does.
    \comm{Here, we are using that for any map \( n \colon U' \to \N \), we have a cover of \( U' \) such that the restriction of \( n \) to each piece is constant.}
  \item The description of when a subobject \( A \tto \code \times \pow(\two^\N) \) satisfies \( \wt \cl_\cap \), \( \wt \cl_\cup \), or \( \wt \cl_\neg \) is similar.
  \item The formula \( \Bo(p) \) defines a subobject \( \Bo \tto \pow(\two^\N) \) so that \( p \colon \wh U \to \pow(\two^\N) \) factors through \( \Bo \) iff \( U \Vdash \Bo(p) \).
    It is the least subobject containing the image of \( i_\basic \) and closed under \( \bigcap \), \( \bigcup \) and \( \neg \).
  \item The formula \( \cBo(q) \) defines an analogous subobject \( \cBo \tto \code \times \pow(\two^\N) \), and \( \uBo(p) \) defines a subobject \( \uBo \tto \pow(\two^\N) \) which is the image of \( \cBo \) under the projection \( \pi_{\pow(\two^\N)} \colon \code \times \pow(\two^\N) \to \pow(\two^\N) \).
  \end{itemize}

  We now proceed with the proof.
  We first prove the converse to plain Borel myopia: that the subobject \( \Bo \) is contained in \( \two^{\two^\N} \tto \pow(\two^\N) \).
  By the definition of \( \Bo \), this amounts to proving that \( \two^{\two^\N} \) is closed under \( \bigcap \), \( \bigcup \), and \( \neg \), and contains the image of \( i_{\basic} \colon \two^{<\N} \to \pow(\two^\N) \).
  But these are all things we have already observed above: that \( \bigcap r,\bigcup r \colon \wh U \to \pow(\two^\N) \) factors through \( \two^{\two^\N} \) whenever \( r \colon \wh U \to \pow(\two^\N)^\N \) factors through \( (\two^{\two^\N})^\N \) (since the Borel subsets are closed under countable union and intersection in \( \Sub(\wh U \times \two^\N) \)); that \( i_{\basic} \colon \two^{<\N} \to \pow(\two^\N) \) factors through \( \two^{\two^\N} \); and the obvious fact that \( \two^{\two^\N} \) is closed in \( \pow(\two^\N) \) under complements.

  We next turn to plain Borel myopia, i.e., that \( \two^{\two^\N} \) is contained in \( \Bo \).
  We need to show that any subobject \( A \tto \pow(\two^\N) \) containing the image of \( i_\basic \) and closed under \( \bigcap \), \( \bigcup \), and \( \neg \) also contains \( \two^{\two^\N} \).
  For each \( U \in \SBor \), consider the set \( \cB_U \) of morphisms \( a \colon U \to \two^{\two^\N} \)---i.e., of Borel subsets \( A \tto \two^\N \times U \)---which factor through \( \Bo \).
  We wish to show that for all \( U \), \( \cB_U \) contains all Borel subsets of \( \two^\N \times U \).
  We will show that \( B \times V \in \cB_U \) for any Borel \( V \subseteq U \) and any basic open \( B \subseteq \two^\N \).
  Since \( \cB_U \) is, by definition of \( \Bo \), closed under countable unions, intersections, and complements, this will finish the proof, since the above sets \( B \times V \) generate the whole Borel sigma algebra.

  First consider the case in which \( V = U \).
  We then have \( B \times V = B \times U \in \cB_U \) for any basic open \( B \), since the morphism \( p \colon U \to \two^\N \) classifying any \( B \times U \tto \two^\N \times U \) factors through \( i_{\basic} \), as observed above.
  In particular, we have \( \two^\N \times U \in \cB_U \), and hence, by closure under complements, that \( \emptyset \in \cB_U \).
  But now for a general \( B \times V \subseteq \two^\N \times U \), we can write it as the disjoint union of \( \two^\N \times V \in \cB_{V} \) and \( \emptyset \times (U \setminus V) \in \cB_{U \setminus V} \).
  If \( a \colon U \to \pow(\two^\N) \) is the morphism classifying \( B \times V \), it follows that the restrictions of \( a \) to \( V \) and \( U \setminus V \) factor through \( \Bo \), and hence that \( a \) does as well, as desired, since \( V \) and \( U \setminus V \) cover \( U. \)

  It remains to show that \( \Bo, \uBo \tto \pow(\two^\N) \) are equivalent subobjects.
  By Proposition~\ref{propn:uniform-implies-plain}, we have the inclusion \( \uBo \le \Bo \).
  To prove the other inclusion, it suffices, by the definition of the \( \Bo \), to show that \( \uBo \) contains the image of \( i_\basic \) and is closed under \( \bigcap \), \( \bigcup \), and \( \neg \).

  The first of these is clear, since \( i_\basic \colon \two^{<\N} \to \pow(\two^\N) \) factors through \( \br{f_{\basic}, i_{\basic}} \colon \two^{<\N} \to \code \times \pow(\two^\N) \), and hence through \( \cBo \tto \code \times \pow(\two^\N) \tox{\pi_{\pow(\two^\N)}} \pow(\two^\N) \).
  Regarding closure under \( f_\neg \), consider \( U \in \SBor \) and a morphism \( p \colon \wh U \to \pow(\two^\N) \) which factors through \( \uBo \), so that there is some epi \( u \colon \wh{U'} \to \wh U \) and \( t \colon \wh{U'} \to \code \) so that \( \br{t, p u} \colon \wh{U'} \to \code \times \pow(\two^\N) \) factors through \( \cBo \).
  We need to show that \( \neg p \colon \wh U \to \pow(\two^\N) \) also factors through \( \uBo \).
  But, first of all, this is true of \( \neg p u \), since \( \br{f_\neg t, \neg p u} \colon \wh{U'} \to \code \times \pow(\two^\N) \) factors through \( \cBo \) by the definition of the latter; and this implies that the same is true of \( \neg p \), since \( u \) is an epimorphism.

  The arguments for closure under \( f_\cap \) and \( f_\cup \) are similar.
  Given a morphism \( p \colon \wh U \to \pow(\two^\N)^\N \)---i.e., a sequence of morphisms \( p_n \colon \wh U \to \pow(\two^\N) \)---such that each \( p_n \) factors through \( \uBo \), we need to show that \( \bigcap p \colon \wh U \to \pow(\two^\N) \) does as well.
  By assumption (and using the countable axiom of choice, externally) we have epimorphisms \( u_n \colon \wh{U_n} \to \wh U \) in \( \SBor \) for each \( n \in \N \) and morphisms \( t_n \colon \wh{U_n} \to \code \) so that each \( \br{t_n, p_n u_n} \colon \wh{U_n} \to \code \times \pow(\two^\N) \) factors through \( \cBo \).
  Now let \( \wh{U'} \) be the fibre product of all the \( \wh{U_n} \) over \( \wh U \), i.e., \( U' = \set{(v_n)_{n \in \N} \in \prod_{n \in \N} U_n \mid u_k(v_k) = u_l(v_l) \FORALL k,l \in \N} \), so that the map \( u = u_n \pi_n \colon \wh{U'} \to \wh U \) is independent of \( n \) and (again by the axiom of choice) an epimorphism.
  We then have morphisms \( t_n' = t_n \pi_n \colon \wh{U'} \to \code \) and \( a_n' = a_n u_n \pi_n \colon \wh{U'} \to \pow(\two^\N) \) such that \( \br{t_n', a_n'} \colon \wh{U'} \to \code \times \pow(\two^\N) \) factors through \( \cBo \) for all \( n \).
  Hence, by definition of \( \cBo \), it follows that \( \br{ \bigcap \br{t_n'}_{n \in \N}, \bigcap \br{a_n'}_{n \in \N} } \colon U' \to \code \times \pow(\two^\N) \) factors through \( \cBo \), and hence that \( \bigcap \br{a_n'}_{n \in \N} \colon U' \to \pow(\two^\N) \) factors through \( \uBo \).
  But since \( \br{a_n'}_{n \in \N} = a u \) and \( u \) is surjective, it follows that \( \bigcap a \) likewise factors through \( \uBo \), as desired.
  The argument for \( \bigcup \) is identical.
\end{proof}

\begin{rmk} \label{rmk:no-strong-myopia}
  Finally, let us now see that, as mentioned above, no topos can satisfy \emph{strong} uniform Borel myopia---i.e., it can never be the case that \emph{every} subset of \( \two^\N \) is uniformly Borel.
  We will explain the argument informally, but it is clear that it only uses intuitionistic higher-order arithmetic, and hence that it is valid in any elementary topos.
  The claim simply amounts to explicitly constructing a specific subset of \( \two^\N \) which is not uniformly Borel, which we do by the usual kind of diagonal argument.

  Note that the composite
  $\cBo \tto \code \times \pow(\two^\N) \to \code$
  is monic
  (this is by induction),
  and fix a monomorphism $\code \tto \two^\N$
  (for instance by fixing a monomorphism
  $\N_3 + \two^{<\N} \tto \two^\N$
  and then using a bijection
  $\N \cong \N \times \N^{<\N}$
  in the exponent).
  Their composite is a monomorphism $f : \cBo \tto \two^\N$.
  Now define
  \[
    u =
    \set{
      (x, y) \in (\two^\N)^2 |
      \exists \br{t, a} \tolon \cBo.\
      (x = f(\br{t, a}) \text{ and }y \in a)
    }.
  \]
  We claim that
  $b = \set{x \in \two^\N : (x, x) \notin u}$
  is not uniformly Borel.
  Suppose it is,
  so that there exists $t \in \code$ with $\br{t, b} \in \cBo$.
  Then setting $x = f(\br{t, b})$,
  by the definition of $u$,
  we have that $(x, x) \in u$ iff $x \in b$
  (we only need to check the second clause since $f$ is monic).
  By definition of $b$,
  however,
  we have $x \in b$ iff $(x, x) \notin U$,
  giving the desired contradiction.
\end{rmk}

\comm{It seems plausible that the Borel topos is the classifying topos for Boolean sigma algebras, but I think this is not so. As I understand, the reason for suspecting it is Chen's observation that \( \SBor \cong \sigma\mathrm{Bool}_\sigma^\op \). In particular, this shows that there is a natural functor \( F \colon \sigma\mathrm{Bool} \to \Fun(\SBor, \Set) \) taking \( B \) to \( X \mapsto \sigma\mathrm{Bool}(\mathcal{B}(X), B) \). Now, if \( \Sh(\SBor) \) were the classifying topos for Boolean \( \sigma \)-algebras, this would also give a functor \( G \colon \sigma\mathrm{Bool} \to \Fun(\SBor, \Set) \), and \emph{presumably}, the idea is that this should be the same as the previous one. But they cannot be the same, since \( G(B) \) will preserve finite limits and take covers to covers, and in particular must take \( \N_{<n} \in \SBor \) to \( \N_{<n} \in \Set \) for all \( n \in \N \), whereas \( F(B) \) does not have this property in general. (For example, when \( B = \mathcal{B}(Y) \) with \( Y \in \SBor \), then \( F(B)(\N_{<n}) \cong \SBor(Y, \N_{<n})) \) by Chen. The possibility does remain that \( \Sh(\SBor) \) is the classifying topos but \( G \) is some \emph{other} functor, but I don't think we have any reason to expect that?}

\appendix

\section{\alttext{\( J_{\SBor} \)}{J\_SBor} is consistently not the canonical topology}\label{sec:necessity-of-ch}

In this appendix, we will prove Proposition~\ref{propn:not_canonical}. That is, we will show that it is consistent with ZFC that $J_\SBor$ is not the canonical topology on $\SBor$. Recall that $J_\SBor$ not being canonical is equivalent to the existence of an uncountable cover $\cU$ of $\Cantor$ by analytic sets with no countable subcover and with the following property: given any Borel map $g\colon X \to \Cantor$ and any function $f\colon X \to Y$ with $X, Y \in \SBor$ , if $f$ is not Borel then there is some $U \in \cU$ such that the restriction of $f$ to $g\I(U)$ is not Borel. Thus we need to show that the existence of such a cover $\cU$ is consistent with ZFC.

In order to prove this, we will describe a model of ZFC in which the cover $\cU$ can be taken to simply be the collection of all Borel subsets of $\Cantor$ of Lebesgue measure $0$ (in fact, in the model we describe we could even take $\cU$ to be just the collection of all $G_\delta$ sets of Lebesgue measure $0$). We will use forcing to construct this model and will use standard terminology and facts about forcing; for a reference on forcing see~\cite{kunen-set-theory-1980}. Before we describe the model, however, we will make a few simplifying observations.

\begin{enumerate}
\item It suffices to consider the case where $X = Y = \Cantor$. That is, it suffices to prove (in the model of ZFC that we construct) that for any Borel map $g\colon \Cantor \to \Cantor$ and any function $f\colon \Cantor \to \Cantor$, if $f$ is not Borel then there is some $U \in \cU$ such that the restriction of $f$ to $g\I(U)$ is not Borel. This follows from the fact that all uncountable standard Borel spaces are Borel isomorphic.

\item It suffices to show that $\cU$ detects non-Borel-ness of sets rather than functions. That is, it suffices to show that for any Borel map $g\colon \Cantor \to \Cantor$ and any set $A \subseteq \Cantor$, if $A$ is not Borel then there is some $U \in \cU$ such that $A\cap g\I(U)$ is not Borel. To see why this is sufficient, suppose that this holds and let $f\colon \Cantor \to \Cantor$ be a function which is not Borel. Then there must be some Borel set $V \subseteq \Cantor$ such that $f\I(V)$ is not Borel. Therefore there must be some $U \in \cU$ such that $f\I(V)\cap g\I(U)$ is not Borel. Since $f\I(V)\cap g\I(U)$ is equal to the preimage of $V$ by $\rstr{f}{g\I(U)}$, this shows that $\rstr{f}{g\I(U)}$ is not Borel.
\end{enumerate}

To summarize, in order to show that it is consistent with ZFC that $J_\SBor$ is not canonical, it is enough to find a model of ZFC such that for every Borel function $g \colon \Cantor \to \Cantor$ and every set $A \subseteq \Cantor$ which is not Borel, there is some Borel set $U$ of Lebesgue measure $0$ such that $A\cap g\I(U)$ is not Borel.

We will now describe a model of ZFC with this property. Let $\bbP$ be the forcing notion which adds $\omega_2$ Cohen reals (i.e.\ $\bbP$ consists of all finite partial functions $\omega_2\times\N \to 2$, ordered by inclusion), let $G$ be $\bbP$-generic over $L$ and let $M = L[G]$. We claim that this model works.\footnote{Actually it doesn't matter that we are forcing over $L$. We could do the same thing over any ground model as long as the number of Cohen reals that we add is larger than the cardinality of the continuum in the ground model.}

In order to prove this, it will help to first establish some notation and recall some standard facts about the model $M$. For each $\alpha < \omega_2$, let $G_{\leq \alpha}$ denote the first $\alpha$ of the Cohen reals added by $G$ and let $G_{> \alpha}$ denote the remaining Cohen reals added by $G$. More precisely, $G_{\leq \alpha} = \rstr{G}{(\alpha + 1)\times \N}$ and $G_{> \alpha} = \rstr{G}{(\alpha + 1, \omega_2)\times \N}$. Let $X_\alpha$ denote $(\Cantor)^{L[G_{\leq \alpha}]}$.

Note that in everything that follows we will always be implicitly working in $M$ unless we explicitly state otherwise.

\begin{lem}
\label{lem:union}
$\Cantor = \bigcup_{\alpha < \omega_2}X_\alpha$.
\end{lem}

\begin{proof}
Let $r \in \Cantor$; we must show that for some $\alpha < \omega_2$, $r \in X_\alpha$. Let $\dot{r}$ be a nice name for $r$ and note that the support of $\dot{r}$ must be countable. Hence there is some $\alpha < \omega_2$ such that the support of $r$ is contained in $\alpha\times \N$. It is straightforward to verify that this implies that $r \in L[G_{\leq \alpha}]$ and hence $r \in X_\alpha$.
\end{proof}

\begin{lem}
\label{lem:measurezero}
For each $\alpha$, $X_\alpha$ has measure $0$.
\end{lem}

\begin{proof}
Note that $G_{> \alpha}$ is Cohen generic over $L[G_{\leq \alpha}]$ and $L[G] = L[G_{\leq \alpha}][G_{> \alpha}]$. It is well known that adding a Cohen generic real turns the ground model reals into a measure $0$ set. Furthermore, the property of a set of reals having measure $0$ is closed upwards under transitive models of ZFC, so in $L[G]$, the set of reals of $L[G_{\leq \alpha}]$, which is precisely the set $X_\alpha$, has measure $0$.
\end{proof}

\begin{lem}
\label{lem:subsetpreimage}
For any Borel function $g\colon \Cantor \to \Cantor$, for all sufficiently large $\alpha$, $X_\alpha \subseteq g\I(X_\alpha)$.
\end{lem}

\begin{proof}
Let $r$ be a Borel code for $g$. By Lemma~\ref{lem:union}, there is some $\beta < \omega_2$ such that $r \in X_\beta$. Now fix $\alpha \geq \beta$ and let $x \in X_\alpha$. Note that $g(x)$ is in $L[x, r]$. Since $x$ and $r$ are both in $L[G_{\leq \alpha}]$, we must have $g(x) \in L[G_{\leq \alpha}]$ and hence $g(x) \in X_\alpha$. Therefore $x \in g\I(X_\alpha)$.
\end{proof}

\begin{thm}
In $M$, $J_\SBor$ is not the canonical topology on $\SBor$.
\end{thm}

\begin{proof}
Working in $M$, let $g\colon \Cantor \to \Cantor$ be a Borel function and let $A \subseteq \Cantor$ be an arbitrary set such that for every Borel set $U \subseteq \Cantor$ of measure $0$, $A \cap g\I(U)$ is Borel. We will show that $A$ is Borel.

We begin with the following two claims, which state that for each $\alpha$, $A \cap X_\alpha$ must be relatively simple. In particular, relative to $X_\alpha$, $A$ is not just Borel but of bounded complexity in the Borel hierarchy (where the bound does not depend on $\alpha$).

\begin{claim}
For each $\alpha < \omega_2$, $A$ is Borel relative to $X_\alpha$. More precisely, there is a Borel set $B$ such that $A\cap X_\alpha = B\cap X_\alpha$.
\end{claim}

\begin{proof}
Fix $\alpha$ and let $\beta > \alpha$ be large enough to apply Lemma~\ref{lem:subsetpreimage}. By Lemma~\ref{lem:measurezero}, $X_\beta$ has measure zero and hence there is some Borel (actually, $G_\delta$) set $U$ of measure $0$ such that $X_\beta \subseteq U$. By our assumption on $A$, $A\cap g\I(U)$ is Borel. We now have
\[
    X_\alpha \subseteq X_\beta \subseteq g\I(X_\beta) \subseteq g\I(U)
\]
where the second inclusion is by Lemma~\ref{lem:subsetpreimage} and our choice of $\beta$. Therefore
\[
    A\cap X_\alpha = A\cap (X_\alpha \cap g\I(U)) =(A \cap g\I(U)) \cap X_\alpha.
\]
Hence $A \cap X_\alpha$ is the intersection of $X_\alpha$ with the Borel set $A\cap g\I(U)$, as desired.
\end{proof}

\begin{claim}
There is some fixed ordinal $\eta < \omega_1$ such that for each $\alpha < \omega_2$, relative to $X_\alpha$, $A$ is of rank at most $\eta$ in the Borel hierarchy. More precisely, for each $\alpha < \omega_2$, there is a $\mathbf{\Sigma}^0_\eta$ set $B$ such that $A\cap X_\alpha = B\cap X_\alpha$.
\end{claim}

\begin{proof}
By the previous claim, for each $\alpha < \omega_2$, we can pick some $\eta_\alpha < \omega_1$ for which there is a $\mathbf{\Sigma}^0_{\eta_\alpha}$ set $B$ such that $A\cap X_\alpha = B\cap X_\alpha$. By the pigeonhole principle, there must be a cofinal set $S\subseteq \omega_2$ on which the value of $\eta_\alpha$ is constant. Let $\eta$ be such a constant value.

Now fix any $\alpha < \omega_2$. Pick some $\beta > \alpha$ which is in $S$. By definition of $\eta$, there is some $\mathbf{\Sigma}^0_\eta$ set $B$ such that $A \cap X_\beta = B\cap X_\beta$. And since $X_\alpha \subseteq X_\beta$, we have that $A\cap X_\alpha = B\cap X_\alpha$ as well.
\end{proof}

For the rest of this proof, let $\eta$ be as in the claim above. Fix a universal $\mathbf{\Sigma}^0_\eta$ set $R \subseteq \Cantor\times\Cantor$. In other words, $R$ is a $\mathbf{\Sigma}^0_\eta$ set such that for each $\mathbf{\Sigma}^0_\eta$ set $B \subseteq \Cantor$, there is some $b \in \Cantor$ such that $B = R_b = \{x \in \Cantor \mid (x, b) \in R\}$. As is standard, we will refer to $b$ as a $\mathbf{\Sigma}^0_\eta$ \defword{code} for $B$. Note that we may assume $R$ is coded in $L$.

Now fix a name $\dot{A}$ for $A$. By modifying $\dot{A}$ (while keeping $\dot{A}[G] = A$), we may assume that $\bbone_\bbP$ forces that for all $\alpha < \omega_2$, $\dot{A}\cap X_\alpha$ is the intersection of $X_\alpha$ with a $\mathbf{\Sigma}^0_\eta$ set. Pick a sequence of nice names $\{\dot{b}_\alpha\}_{\alpha < \omega_2}$ giving $\mathbf{\Sigma}^0_\eta$ codes for these sets. More precisely, pick a sequence such that for all $\alpha$, $\bbone_\bbP$ forces that $\dot{b}_\alpha$ is a $\mathbf{\Sigma}^0_\eta$ code for a $\mathbf{\Sigma}^0_\eta$ set $B$ such that $\dot{A} \cap X_\alpha = B \cap X_\alpha$.

Our goal is to somehow amalgamate these names into a single Borel code for the set $A$. We will begin by restricting to a cofinal subset of $\omega_2$ on which the names $\dot{b}_\alpha$ behave particularly nicely. In order to do so, it helps to establish some notation.

For each condition $p$ in $\bbP$, the \defword{support} of $p$, denoted $\supp(p)$, is the set of ordinals $\alpha < \omega_2$ such that for some $n \in \N$, $(\alpha, n)$ is contained in the domain of $p$. For any name $\dot{x}$, the \defword{support} of $x$, denoted $\supp(\dot{x})$,  is the union of the supports of conditions which appear in $\dot{x}$. Given a condition $p \in \bbP$ and a map $\phi \colon \supp(p) \to \omega_2$, let $\phi(p)$ denote the condition defined by $\phi(p)(\alpha, n) = p(\phi(\alpha), n)$. Similarly, for a name $\dot{x}$ and a map $\phi \colon \supp(\dot{x}) \to \omega_2$, let $\phi(\dot{x})$ denote the name obtained by replacing each condition $p$ in $\dot{x}$ by $\phi(p)$.

Note that the support of each name $\dot{b}_\alpha$ is countable (since it is a nice name for a real). Hence by applying the $\Delta$-system lemma to the supports of the $\dot{b_\alpha}$, we obtain a countable set $C \subseteq \omega_2$ and a cofinal set $S_1 \subseteq \omega_2$ such that for all $\alpha \neq \beta$ in $S_1$, $\supp(\dot{b}_\alpha) \cap \supp(\dot{b}_\beta) = C$ (i.e.\ $C$ is the stem of the $\Delta$-system formed by $\{\supp(\dot{b}_\alpha) \mid \alpha \in S_1\}$).

For each $\alpha \in S_1$, let $D_\alpha$ denote the set $\supp(\dot{b}_\alpha)\setminus C$ (i.e. the part of $\supp(\dot{b}_\alpha)$ that is outside the stem $C$). It is easy to see that we may refine $S_1$ to a cofinal set $S_2 \subseteq \omega_2$ such that if $\alpha < \beta$ are in $S_2$ then every element of $D_\alpha$ is below every element of $D_\beta$.

Now note that each $D_\alpha$ is countable and hence has order type less than $\omega_1$. By the pigeonhole principle, we may restrict $S_2$ to a cofinal subset $S_3$ so that for all $\alpha, \beta \in S_3$, $D_\alpha$ and $D_\beta$ have the same order type.

Let $\gamma_0 = \sup\{\alpha + 1 \mid \alpha \in C\}$ and let $\gamma_1$ denote the unique order type of $D_\alpha$ for $\alpha \in S_3$. Let $\gamma = \gamma_0 + \gamma_1$. For each $\alpha \in S_3$, let $\phi_\alpha \colon \supp(\dot{b}_\alpha) \to \gamma$ be the map given by preserving $C$ and collapsing $D_\alpha$ to its order type. More precisely, for each $\beta \in C$, $\phi_\alpha(\beta) = \beta$ and for each $\beta \in D_\alpha$, $\phi_\alpha(\beta) = \gamma_0 + \beta'$ where $\beta'$ is the image of $\beta$ when $D_\alpha$ is collapsed to its order type, $\gamma_1$.

Note that for any $\alpha \in S_3$, $\phi_\alpha(\dot{b}_\alpha)$ is a countable name with support $\gamma < \omega_1$ and hence can be coded by a real. So as $\alpha$ ranges over $S_3$ there are at most $|\two^\N|^L = \omega_1$ possible values for $\phi_\alpha(\dot{b}_\alpha)$. Therefore by applying the pigeonhole principle again, we may restrict to a cofinal set $S_4 \subseteq S_3$ such that $\phi_\alpha(\dot{b}_\alpha)$ is constant on $S_4$. Let $\dot{b}$ denote this constant value.

The idea is that $\dot{b}$ more or less gives us a $\mathbf{\Sigma}^0_\eta$ code for $A$. However, it takes a bit of work to make this precise.

Let $\bbP_{\gamma}$ denote the forcing that adds $\gamma$-many Cohen reals and let $\bbP_{\gamma_1}$ denote the forcing that adds $\gamma_1$-many Cohen reals. Note that $\dot{b}$ can be thought of as a nice $\bbP_\gamma$-name for an element of $\Cantor$. In particular, $\dot{b}$ is a set of pairs of the form $(n, p)$ where $n \in \N$ and $p$ is a $\bbP_\gamma$-condition. Let $\dot{c}$ denote the name given by excluding from $\dot{b}$ all pairs $(n, p)$ such that $p$ is not compatible with $G_{\gamma_0}$. Note that $\dot{c}$ can be thought of as a nice $\bbP_{\gamma_1}$-name for an element of $\Cantor$.

Let $B$ be the set of $x \in \Cantor$ such that $\bbone_{\bbP_{\gamma_1}}$ forces that $x$ is in the $\mathbf{\Sigma}^0_\eta$ set coded by $\dot{c}$. We claim that $B$ is a Borel set and that it is equal to $A$.

\begin{claim}
The set $B$ is Borel.
\end{claim}

\begin{proof}
Note that a generic filter for $\bbP_{\gamma_1}$ can be coded by a real. Furthermore, Cohen forcing admits continuous reading of names so there is a continuous function $F$ such that for any $\bbP_{\gamma_1}$-generic $h$, $F(h) = \dot{c}[h]$. Hence the question of whether $x$ is in $B$ is equivalent to the question of whether it is forced that $(x, F(h)) \in R$. Since $R$ is a $\mathbf{\Sigma}^0_\eta$ set coded in $M$, the standard analysis of Cohen forcing shows that the property
\[
    \bbone_{\bbP_{\gamma_1}} \Vdash (x, F(h)) \in R
\]
is Borel.
\end{proof}

\begin{claim}
$A = B$.
\end{claim}

\begin{proof}
Fix $x \in \Cantor$. We will show that if $x \in A$ then $x \in B$. The proof that if $x \notin A$ then $x \notin B$ is very similar and so we omit it.

Assume that $x \in A$. By Lemma~\ref{lem:union} we can pick $\alpha < \omega_2$ such that $x \in X_\alpha$. As we have observed above, $L[G] = L[G_{\leq \alpha}][G_{> \alpha}]$. Since $x \in A$ holds in $L[G_{\leq \alpha}][G_{> \alpha}]$, $\check{x} \in \dot{A}$ must be forced over $L[G_{\leq \alpha}]$ by some condition in $G_{> \alpha}$. Let $p$ be such a condition.

Note that $\supp(p)$ is finite. Hence we may pick some $\alpha' < \omega_2$ larger than both $\alpha$, $\gamma_0$ and every element of $\supp(p)$. Let $\beta \in S_4$ be large enough that every element of $D_\beta$ is larger than both $\alpha'$ (recall that $S_4$ is cofinal in $\omega_2$ and all the $D_\beta$'s are countable and disjoint). Let $\bbP_{\leq \alpha'}$ and $\bbP_{> \alpha'}$ denote the obvious forcing notions so that $G_{\leq \alpha'}$ is generic for $\bbP_{\leq \alpha'}$ and $G_{> \alpha'}$ is generic for $\bbP_{> \alpha'}$ over $L[G_{\leq \alpha'}]$.

Since $\alpha'$ is larger than every element of $\supp(p)$ and since $G_{\leq \alpha'}$ extends $p$ (because $p$ is in $G_{> \alpha}$ and $G_{\leq \alpha'}$ agrees with $G_{> \alpha}$ on $\supp(p)$), we must have that in $L[G_{\leq \alpha'}]$, $\bbone_{\bbP_{> \alpha'}}$ forces that $\check{x} \in \dot{A}$. Hence by our choice of $\dot{b}_{\beta}$ it must also be forced by $\bbone_{\bbP_{\alpha'}}$ that $\check{x}$ is in the $\mathbf{\Sigma}^0_\eta$ set coded by $\dot{b}_\beta$.

Note that $\bbP_{> \alpha'}$ is isomorphic to $\bbP_{> \alpha'}\times \bbP_{\gamma_1}$ by mapping $D_\beta$ to $\gamma_1$. Moreover, this mapping sends $\dot{b}_\beta$ to $\dot{c}$. Therefore we have that in $L[G_{\leq \alpha'}]$, it is forced by $\bbone_{\bbP_{> \alpha'}\times\bbP_{\gamma_1}}$ that $\check{x}$ is in the $\mathbf{\Sigma}^0_\eta$ set coded by $\dot{c}$. Since $G_{> \alpha'}$ is $\bbP_{> \alpha'}$-generic, this implies that in $L[G] = L[G_{\leq \alpha'}][G_{> \alpha'}]$ it is forced by $\bbone_{\bbP_{\gamma_1}}$ that $\check{x}$ is in the $\mathbf{\Sigma}^0_\eta$ set coded by $\dot{c}$.
\end{proof}

Clearly the two claims above are enough to finish the proof.
\end{proof}

\section{Proof of Theorem~\ref{thm:reg-completion}} \label{sec:reg-completion-proof}
We need to show that \( \SAn \) is the regular completion of \( \SBor \), the regularly countably complete, countably extensive completion of \( \SBor \), and the completion of \( \SBor \) in all of the intermediate senses explained in \S\ref{subsec:analytic}.

To begin with, that the inclusion \( \SBor \hto \SAn \) is regular, countably complete, and countably extensive follows from Proposition~\ref{propn:san-is-good}.

Now let \( \cC \) be a regular category and \( F \colon \SBor \to \cC \) a regular functor.
We must first show this extends to a regular functor \( \wt F \colon \SAn \to \cC \).
For each \( A \in \SAn \), fix an inclusion \( i_A \colon A \hto Y_A \in \SBor \) and a surjective morphism \( q_A \colon X_A \to A \) with \( X_A \in \SBor \) (taking \( i_A = q_A = \id_A \) if \( A \) is Borel), choose a regular epi-mono factorization (i.e., an image) \( F X_A \tox{\tilde q_A} \wt A \tox{\tilde i_A} F Y_A  \) in \( \cC \) of \( F(i_A \circ q_A) \colon F X_A \to F Y_A \), and set \( \wt F A = \wt A \), again taking \( \wt A = F A \) and \( \tilde q_A = \tildei_A = \id_{F A} \) for \( A \) Borel.\footnote{
  The use of choice is a common occurrence when proving universal properties free constructions in category theory, but is inessential, in the sense that it can be avoided by the use of \emph{anafunctors} in the sense of \cite{makkai-avoiding-choice}.
}
Next, for each morphism \( f \colon A \to B \), we claim that there is a unique morphism \( \tilde f \colon \wt F A \to \wt F B \) with \( \tildei_B \tilde f q_A = F(i_B f q_A) \colon F X_A \to F Y_B \), and we will then define \( \wt F f = \tilde f \).
The uniqueness of \( \tilde f \) is immediate since \( \tildei_B \) is monic \( \tilde q_A \) is epic.
\[
  \begin{tikzcd}
    X_A \ar[r, "q_A", ->>] & A \ar[r, "i_A", >->] \ar[d, "f"] & Y_A \\
    X_B \ar[r, "q_B", ->>] & B \ar[r, "i_B", >->] & Y_B
  \end{tikzcd}
  \qquad
  \begin{tikzcd}
    F X_A \ar[r, "\tilde q_A", ->>] & \wt F A \ar[r, "\tildei_A", >->] \ar[d, "\tilde f", dashed] & F Y_A \\
    F X_B \ar[r, "\tilde q_B", ->>] & \wt F B \ar[r, "\tildei_B", >->] & F Y_B
  \end{tikzcd}
\]
To prove that \( F(i_B f q_A) \) factors through \( \wt F B \), we form the following pullback in \( \SBor \), noting that \( \pi_{X_A} \) is a regular epi by Proposition~\ref{propn:san-is-good}, and hence we get a corresponding pullback square and regular epi in \( \cC \).
\begin{equation} \label{eq:san-sbor-pbs}
  \begin{tikzcd}[column sep=40pt]
    P \pb[rrd] \ar[r, "\pi_{X_A}", ->>] \ar[rd, "\pi_{X_B}"'] & X_A \ar[r, "i_A q_A"] \ar[rd, "i_B f q_A"] & Y_A \\
    & X_B \ar[r, "i_B q_B"] & Y_B
  \end{tikzcd}
  \qquad
  \begin{tikzcd}[column sep=40pt]
    F P \pb[rrd] \ar[r, "F \pi_{X_A}", ->>] \ar[rd, "F \pi_{X_B}"'] & F X_A \ar[r, "\tildei_A \tilde q_A"] \ar[rd, "F(i_B f q_A)"] & F Y_A \\
    & F X_B \ar[r, "\tildei_B \tilde q_B"] & F Y_B.
  \end{tikzcd}
\end{equation}
Since \( F(i_B f q_A) \pi_{X_A} \) factors through the mono \( \tildei_B \) and \( F \pi_{X_A} \) is a regular epi, it follows that \( F(i_B f q_A) \) also factors through \( \tildei_B \), as desired.

Next, to see that \(  F(i_B f q_A) \) factors through \( \wt F A \), it suffices, since \( \tilde q_A \) is a regular epi, to show that \( F(i_B f q_A) \pi_0 = F(i_B f q_A) \pi_1 \colon F X_A \times_{\wt F A} F X_A \to F Y_B \).
Since \( \tildei_A \) is a mono and \( F \) preserves finite limits, the canonical morphisms
\[
  F(X_A \times_{Y_A} X_A) \to F X_A \times_{F Y_A} F X_A \ot F X_A \times_{\wt F A} F X_A
\]
are isos, so we can work with the pullback \( F(X_A \times_{Y_A} X_A) \) instead.
But since \( i_B f q_A \) factors through the image \( A \) of \( i_A q_A \), it follows that \( i_B f q_A \pi_0 = i_B f q_A \pi_1 \colon X_A \times_{Y_A} X_A \to Y_B \) in \( \SBor \), and hence that \( F(i_B f q_A) F(\pi_0) = F(i_B f q_A) F(\pi_1) \colon F(X_A \times_{Y_A} X_A) \to F(Y_B) \), as desired.

We have established the well-definedness of \( \tilde f = \wt F f \); we next need to show that \( \wt F \) is a functor.
That it preserves identity morphisms is clear.
To see that it preserves composition, consider morphisms \( A \tox{f} B \tox{g} C \) in \( \SAn \), giving diagrams
\[
  \begin{tikzcd}
    P \ar[r, "\pi_{X_A}", ->>] \ar[rd, "\pi_{X_B}"'] & X_A \ar[r, "q_A", ->>] & A \ar[r, "i_A", >->] \ar[d, "f"] & Y_A \\
    & X_B \ar[r, "q_B", ->>] & B \ar[r, "i_B", >->] \ar[d, "g"] & Y_B \\
    & X_C \ar[r, "q_C", ->>] & C \ar[r, "i_C", >->] & Y_C
  \end{tikzcd}
  \qquad
  \begin{tikzcd}
    F P \ar[r, "F \pi_{X_A}", ->>] \ar[rd, "F \pi_{X_B}"'] &[5pt] F X_A \ar[r, "\tilde q_A", ->>] & \wt F A \ar[r, "\tildei_A", >->] \ar[d, "\wt F f", dashed] & F Y_A \\
    & F X_B \ar[r, "\tilde q_B", ->>] & \wt F B \ar[r, "\tildei_B", >->] \ar[d, "\wt F g", dashed] & F Y_B \\
    & F X_C \ar[r, "\tilde q_C", ->>] & \wt F C \ar[r, "\tildei_C", >->] & F Y_C,
  \end{tikzcd}
\]
where we have formed the same pullbacks as \eqref{eq:san-sbor-pbs} in above.
To prove that \( \wt F g \cdot \wt F g = \wt F(g f) \), it suffices to prove this after composing with the epi \( \tilde q_A \cdot F \pi_{X_A} \) and the mono \( \tildei_C \).
First, we have \( \tildei_B \cdot \wt F f \cdot \tilde q_A \cdot F \pi_{X_A} = \tildei_B \cdot \tilde q_B \cdot F \pi_{X_B} \) by definition of \( P \), and hence  \( \wt F f \cdot \tilde q_A \cdot F \pi_{X_A} = \tilde q_B \cdot F \pi_{X_B} \) since \( \tildei_B \) is mono.
Using the definition of \( \wt F \), we then have:
\[
  \tildei_C \cdot \wt F g \cdot \wt F f \cdot \tilde q_A \cdot F \pi_{X_A}
  =
  \tildei_C \cdot \wt F g \cdot \tilde q_B \cdot F \pi_{X_B}
  =
  F(i_C g q_B) \cdot F \pi_{X_B}
  =
  F(i_C g q_B \pi_{X_B})
\]
and
\[
  \tildei_C \cdot \wt F(g f) \cdot \tilde q_A \cdot F \pi_{X_A} =
  F(i_C g f q_A) \cdot F \pi_{X_A} =
  F(i_C g f q_A \pi_{X_A}).
\]
It thus remains to prove that \( i_C g q_B \pi_{X_B} = i_C g f q_A \pi_{X_A} \) in \( \SBor \), or equivalently in \( \SAn \).
But this follows since \( q_B \pi_{X_B} = f q_A \pi_{X_A} \) in \( \SAn \), because \( i_B q_B \pi_{X_B} = i_B f q_A \pi_{X_A} \) in \( \SBor \), hence also in \( \SAn \), and \( i_B \) is monic.

This finishes the definition of the functor \( \wt F \).
We next need to see that \( \wt F \), thus defined, is regular.
First, it is immediate that for each morphism \( q_A \colon X_A \to A \) which was chosen in the definition of \( F \), we have \( \wt F q_A = \tilde q_A \), which is by definition a regular epi.
Then, any epi \( f \colon X \to Y \) in \( \SAn \) with \( X \) Borel is isomorphic to such a \( q_A \), hence \( \wt F f \) is isomorphic to \( \tilde q_A \) and hence regular epi as well.
Finally, for a general \( g \colon A \to B \) in \( \SAn \), we can choose an epi \( f \colon X \to A \) with \( X \) Borel; then \( \wt F g \cdot \wt F f = \wt F(g f) \colon X \to B \) is a regular epi by what we just showed, and hence so is \( \wt F g \).

It remains to see that \( \wt F \) preserves finite limits.
That it preserves the terminal object is clear.
We need to see that it preserves pullbacks.
\[
  \begin{tikzcd}
    X' \pb \ar[r, ->>] \ar[d] & A' \pb \ar[r] \ar[d]& B' \pb \ar[r, equals] \ar[d] & B' \ar[d] \\
    X \ar[r, ->>] & A \ar[r] & B \ar[r, >->] & Y
  \end{tikzcd}
\]
Consider a pullback square in \( \SAn \) with vertices \( A,B,A',B' \) as shown above.
We may reduce to the case in which \( B \in \SBor \) by pasting on the square shown above to the right with \( Y \in \SBor \), and using the usual 2-of-3 rule for pullback squares.
We may then further reduce to the case in which \( A \in \SBor \) by pasting on the square shown above to the left and using that, in a regular category, the ``unusual'' 2-of-3 rule (i.e., left and outer squares are pullbacks imply right square is a pullback) holds if the horizontal arrows in the left square are regular epis.\comm{This take a bit of doing, but is fairly straightforward, and seems to be well-known (e.g., is discussed on stackexchange).}
By the same device, we may further reduce to the case in which \( B' \), and hence also \( A' \) is in \( \SBor \).
And then we are done since \( \wt F \big|_{\SBor} = F \) preserves finite limits by assumption.

We have shown that any regular \( F \colon \SBor \to \cC \) extends to a regular \( \wt F \colon \SAn \to \cC \).
Next, we must show that any natural transformation \( \alpha \colon F \to F' \) between functors \( F,F \colon \SBor \to \SAn \) extends to a unique natural transformation \( \tilde \alpha \colon \wt F \to \wt F' \) between given regular extensions \( \wt F, \wt F' \colon \SAn \to \cC \).
Given \( A \in \SAn \) choose an epi and mono \( X \tox{q_A} A \tox{i_A} Y \) with \( X,Y \in \SBor \), as above, and consider the diagram
\[
  \begin{tikzcd}[column sep=40pt]
    F X_A \ar[r, "F q_A", ->>] \ar[d, "\alpha_{X_A}"] & \wt F A \ar[r, "F i_A", >->] \ar[d, "\tilde \alpha_A", dashed] & \wt F Y_A \ar[d, "\alpha_{Y_A}"] \\
    F' X_A \ar[r, "F' q_A", ->>] & \wt F' A \ar[r, "F' i_A", >->] & F' Y_A.
  \end{tikzcd}
\]
We then observe that the outer square commutes by naturality of \( \alpha \), and define \( \tilde \alpha_A \) to be the unique dashed morphism making the whole diagram commute.
Moreover, this is the only possible choice by the required naturality, either with respect to \( q_A \) or with respect to \( i_A \).
That \( \tilde \alpha \) is natural with respect to any given morphism \( f \colon A \to B \) in \( \SAn \) follows by a similar argument using the naturality of \( \alpha \) with respect to \( i_B f q_A \colon X_A \to Y_B \).

It remains to see that, if \( \cC \) additionally has any of regularly countable products, finite stable disjoint coproducts, countable stable disjoint coproducts, finite stable joins, or countable stable joins and, in each case, if \( F \) preserves these, then so does \( \wt F \).

For regularly countable products: give \( \prod_{n \in \N} A_n \in \SAn \), choose \( X_{A_n} \too A_n \tto Y_{A_n} \) for each \( n \); this gives induced maps \( \prod_n X_{A_n} \too \prod_n A_n \tto \prod_n Y_{A_n} \) and hence a diagram in
\[
  \begin{tikzcd}
    F\prod_{n \in \N} X_n \ar[r, ->>] \ar[d, "\sim" sloped] &
    \wt F\prod_{n \in \N} A_n \ar[r, >->] \ar[d] &
    F\prod_{n \in \N} Y_n \ar[d, "\sim" sloped]
    \\
    \prod_{n \in \N} F X_n \ar[r, ->>] &
    \prod_{n \in \N} \wt F A_n \ar[r, >->] &
    \prod_{n \in \N} F Y_n
  \end{tikzcd}
\]
in \( \cC \), where we use the regularly countable completeness of \( \cC \) to conclude that the lower left horizontal morphism is a regular epi.
Since the two outer vertical morphisms are isos, the central one is as well, as desired.
The argument for stable disjoint coproducts is similar, using that both regular epis and monos are closed under arbitrary stable disjoint coproducts. \comm{This takes a bit of doing, but is fairly straightforward.}

Finally, we consider stable joins; the argument in the finite and countable cases are identical.
Let \( i_n \colon A_n \tto A \) be a (finite or countable) family of monomorphisms in \( \SAn \) whose join in \( \Sub(A) \) is \( A \)---i.e., a covering of \( A \) by analytic subsets.
We must show that the family \( F i_n \colon F A_n \tto F A \) likewise has join \( F A \).
Choose surjections \( u_n \colon X_n \too A_n \) with \( X_n \in \SBor \), thus giving a surjection \( u = [i_n u_n]_n \colon \coprod_n X_n \too A \).
Now suppose we have some \( i_B \colon B \tto F A \) through which each \( F A_n \) factors, and consider the following diagram for each \( n \), where we have pulled back \( B \) along \( F u \):
\[
  \begin{tikzcd}
    &[-10pt] F X_n \ar[rr, "F u_n", ->>] \ar[ld, dashed, >->] &&[-10pt] F A_n \ar[dd, >->, "F i_n"] \ar[ld, >->] \\[-10pt]
    (F u)^* B \ar[rr] \ar[rd, >->, "i_{(F u)^* B}"'] \pb[rrrd] && B \ar[rd, >->, "i_B"] \\[-10pt]
    & F(\coprod_n X_n) \ar[from=uu, crossing over, >->] \ar[rr, "F u", ->>] & {} & FA.
  \end{tikzcd}
\]
Then each \( F X_n \) factors through \( (F u)^* B \) as shown by the universal property of the latter and commutativity of the displayed square.
Since \( \coprod_n X_n \) is the join of the \( X_n \), the same holds for the \( F X_n \tto F(\coprod_n X_n) \) by assumption, and hence, since these all factor through \( (F u)^* B \), we conclude that \( i_{(F u)^* B} \) is an iso, as indicated.
Since \( F u \) is a regular epi, it follows that \( i_B \) is also an iso, as desired.

\printbibliography
\end{document}